\documentclass[paper=a4,9pt]{article}
\usepackage{blindtext}
\title{The free Eggert-operad and an operadic description of groups}
\date{}
\author{Roman Haak\\ University of Stuttgart}
\usepackage{fancyvrb}
\usepackage[utf8]{inputenc}
\usepackage[OT2,T1]{fontenc}
\DeclareSymbolFont{cyrletters}{OT2}{wncyr}{m}{n}
\DeclareMathSymbol{\sha}{\mathalpha}{cyrletters}{"58}
\DeclareSymbolFont{Shuffle}{U}{shuffle}{m}{n}
\DeclareMathSymbol\shuffle{\mathbin}{Shuffle}{"001}
\usepackage{lmodern}
\usepackage{amssymb}
\usepackage{amsmath}
\usepackage{amsthm}
\usepackage{thmtools}
\theoremstyle{definition}
\usepackage[pdftex]{graphicx}
\usepackage{tikz}
\usepackage{tikz-cd}
\usetikzlibrary{
  knots,
  hobby,
  decorations.pathreplacing,
  shapes.geometric,
  calc
}
\tikzset{
  knot diagram/every strand/.append style={
    ultra thick,
    red
  },
  show curve controls/.style={
    postaction=decorate,
    decoration={show path construction,
      curveto code={
        \draw [blue, dashed]
        (\tikzinputsegmentfirst) -- (\tikzinputsegmentsupporta)
        node [at end, draw, solid, red, inner sep=2pt]{};
        \draw [blue, dashed]
        (\tikzinputsegmentsupportb) -- (\tikzinputsegmentlast)
        node [at start, draw, solid, red, inner sep=2pt]{}
        node [at end, fill, blue, ellipse, inner sep=2pt]{}
        ;
      }
    }
  },
  show curve endpoints/.style={
    postaction=decorate,
    decoration={show path construction,
      curveto code={
        \node [fill, blue, ellipse, inner sep=2pt] at (\tikzinputsegmentlast) {}
        ;
      }
    }
  }
}
\usepackage{textcomp}
\usepackage{makeidx}\makeindex
\usepackage[]{geometry}
\usepackage{pgfplots}
\usepackage{filecontents}
\usepackage{url}
\usepackage{etoolbox}
\appto\UrlBreaks{\do\a\do\b\do\c\do\d\do\e\do\f\do\g\do\h\do\i\do\j
\do\k\do\l\do\m\do\n\do\o\do\p\do\q\do\r\do\s\do\t\do\u\do\v\do\w
\do\x\do\y\do\z}
\usepackage{hyperref}
\usepackage{relsize}
\usepackage{stackengine}
\stackMath
\usetikzlibrary{decorations.markings}
\newtheorem{defi}{Definition}[section]
\newtheorem{satz}[defi]{Theorem}
\newtheorem{bei}[defi]{Example}
\newtheorem{lem}[defi]{Lemma}
\newtheorem{bem}[defi]{Remark}
\newtheorem{kor}[defi]{Corollary}
\newtheorem{eri}[defi]{Reminder}
\newtheorem{propo}[defi]{Proposition}
\newtheorem{defilem}[defi]{Definiton/Lemma}
\newtheorem{nota}[defi]{Notation}
\newcommand{\bfcite}[2]{{\bf\cite[\rm #2]{#1}}}
\newcommand{\bfcit}[1]{{\bf\cite{#1}}}

\newcommand{\mapo}{\mathrm{Map}_0}
\newcommand{\mapoop}{\mathrm{Map}_0^{\mathrm{op}}}

\newcommand{\eno}{\mathrm{End}_0}

\newcommand{\Eno}{\mathrm{END}_0}

\newcommand{\Groupfp}{\mathrm{GROUP}_0^{\mathrm{fp}}}
\newcommand{\Groupfppre}{\mathrm{GROUP}_0^{\mathrm{fp,pre}}}

\newcommand{\Wordso}{\mathrm{WORDS}_0(X)}

\newcommand{\Words}{\mathrm{WORDS}(X)}

\newcommand{\Freeox}{\mathrm{FREE}_0(X)}

\newcommand{\Freeoprex}{\mathrm{FREE}_0^{\mathrm{pre}}(X)}

\newcommand{\Freeog}{\mathrm{FREE}_0(X_{\operatorname{G}})}
\newcommand{\Freeopreg}{\mathrm{FREE}^{\mathrm{pre}}_0(X_{\operatorname{G}})}
\newcommand{\Z}{\mathbb{Z}_{\geq 0}}

\newcommand{\perm}{\operatorname{S}}

\newcommand{\id}{\operatorname{id}}

\newcommand{\proppo}{\mathcal{P}_0}

\newcommand{\oppo}{(\proppo,\mathfrak{p}_0)}

\newcommand{\qproppo}{\mathcal{Q}_0}
\newcommand{\qpropp}{\mathcal{Q}}
\newcommand{\qoppo}{(\qproppo,\mathfrak{q}_0)}

\newcommand{\rproppo}{\mathcal{R}_0}

\newcommand{\xy}{\big<X\,|\,Y\big>}

\newcommand{\xgyg}{\big<X_{\operatorname{G}}\,|\,Y_{\operatorname{G}}\big>}

\newcommand{{\idealyGamma}}{\big< Y\mapGamma \big>}

\newcommand{\generated}{{}_{\mathrm{op}}}

\newcommand{\urho}{\uprho}

\newcommand{\fzero}{f^{(0)}}
\newcommand{\ftwo}{f^{(2)}}

\newcommand{\biifree}{\mathsf{fr}} % \mathbf{f}\quad \mathrm{f} \quad \mathsf{f} \quad \mathtt{f}
\newcommand{\biipropp}{\mathsf{p}} % \mathbf{f}\quad \mathrm{f} \quad \mathsf{f} \quad \mathtt{f}
\newcommand{\mapGamma}{\Gamma}

\newcommand{\blue}{\textcolor{blue}}

\newcommand{\fp}{^\mathrm{fp}}
\newcommand{\op}{^\mathrm{op}}

\newcommand{\catob}{\mathrm{Ob}}

\newcommand{\funcfacfunc}{\urho}

\newcommand{\so}{\mathrm{s}}
\newcommand{\ta}{\mathrm{t}}

\newcommand{\IV}{\text{IH}}
\newcommand{\cf}{cf.\,}
\newcommand{\cfeg}{cf.\hspace{-0.5ex} e.g.}
\newcommand{\Cf}{Cf.\,}
\newcommand{\nV}{\text{by ass.}}
\newcommand{\diag}{\mathrm{diag}}
\newcommand{\konst}{\mathrm{const}}
\usepackage{stmaryrd}%für umgedrehtes \mapsto, Befehl heißt \mapsfrom
\newcommand{\se}{\mathsf{s}}
\newcommand{\he}{\mathsf{h}}
\newcommand{\ie}{\operatorname{i}\hspace{-0.5ex}.\hspace{-0.5ex}\operatorname{e}\hspace{-0.5ex}.\hspace{-0.2ex}}

\newcommand{\eg}{\operatorname{e}\hspace{-0.5ex}.\hspace{-0.5ex}\operatorname{g}\hspace{-0.5ex}.\hspace{-0.2ex}}
\newcommand{\vsp}{\vspace{1ex}}

\newcommand{\bq}{\begin{quote}\begin{footnotesize}}
\newcommand{\eq}{

                 \end{footnotesize}\end{quote}} % die Leerzeiche ist Absicht, da sonst Formatierungsprobleme
\newcommand{\bdefi}{\begin{defi}}
\newcommand{\edefi}{

                 \end{defi}} % die Leerzeiche ist Absicht, da sonst Formatierungsprobleme
                 \newcommand{\blem}{\begin{lem}}
\newcommand{\elem}{

                 \end{lem}} % die Leerzeiche ist Absicht, da sonst Formatierungsprobleme
                                  \newcommand{\bpro}{\begin{proof}}
\newcommand{\epro}{\end{proof}} % die Leerzeiche ist Absicht, da sonst Formatierungsprobleme
                 \newcommand{\bsatz}{\begin{satz}}
\newcommand{\esatz}{

                 \end{satz}} % die Leerzeiche ist Absicht, da sonst Formatierungsprobleme
                 \newcommand{\bbem}{\begin{bem}}
\newcommand{\ebem}{

                 \end{bem}} % die Leerzeiche ist Absicht, da sonst Formatierungsprobleme
                 \newcommand{\bkor}{\begin{kor}}
\newcommand{\ekor}{

                 \end{kor}} % die Leerzeiche ist Absicht, da sonst Formatierungsprobleme
                 \newcommand{\bbei}{\begin{bei}}
\newcommand{\ebei}{

                 \end{bei}} % die Leerzeiche ist Absicht, da sonst Formatierungsprobleme
                 \newcommand{\beri}{\begin{eri}}
\newcommand{\eeri}{

                 \end{eri}} % die Leerzeiche ist Absicht, da sonst Formatierungsprobleme
                \newcommand{\bpropo}{\begin{propo}}
\newcommand{\epropo}{

                 \end{propo}} % die Leerzeiche ist Absicht, da sonst Formatierungsprobleme
                 \newcommand{\bdefilem}{\begin{defilem}}
\newcommand{\edefilem}{

                 \end{defilem}} % die Leerzeiche ist Absicht, da sonst Formatierungsprobleme
                 \newcommand{\bnota}{\begin{nota}}
\newcommand{\enota}{

                 \end{nota}} % die Leerzeiche ist Absicht, da sonst Formatierungsprobleme
\usepackage{upgreek}
\VerbatimFootnotes
\usetikzlibrary{arrows,chains,positioning,scopes,matrix}
\makeatletter
\tikzset{join/.code=\tikzset{after node path={%
\ifx\tikzchainprevious\pgfutil@empty\else(\tikzchainprevious)%
edge[every join]#1(\tikzchaincurrent)\fi}}}
\makeatother
\tikzset{>=stealth',every on chain/.append style={join},
         every join/.style={->}}
\tikzstyle{labeled}=[execute at begin node=$\scriptstyle,
   execute at end node=$]  
\DeclareMathOperator{\btimes}{\pmb{\pmb{\times}}}
\begin{document}

\maketitle

\bq
An Eggert-operad is a variant of Mac Lane's notion of a PROP, for which not only bijective maps, but all maps between standard finite sets, are part of the structure. We construct the free Eggert-operad and prove the universal property it satisfies. Further we use the free Eggert-operad to define the Eggert-operad $\Groupfp$ via generators and relations and show that a $\Groupfp$-algebra is a group and vice versa. At the end we prove operadic versions of some elementary facts about groups.
\eq

\begin{footnotesize} \renewcommand{\baselinestretch}{0.9} \parskip0.0ex \tableofcontents  \parskip1.2ex \renewcommand{\baselinestretch}{1.0} \end{footnotesize}

\section{Introduction}

In the sixties, Mac Lane introduced PROPs to describe algebraic structures such as $\eg$ monoids and algebras from a categorical standpoint \bfcite{MacLane}{p.\ 97}, \bfcit{MacLanePACTs}. A PROP is given by a strict monoidal category $\rproppo$ with $(\Z,+)$ as monoid of objects such that the symmetric group $\perm_m$ is a subgroup of $\rproppo(m,m)$ for $m\in \Z$. Roughly speaking, an element in $\rproppo(m,n)$, where $m,n\in \Z$, gives rise to an operation with $m$ ``inputs'' and $n$ ``outputs''. Boardman and Vogt used PROPs to study infinite loop spaces \bfcite{BV}{p.\ 1119}. 

In the seventies May introduced operads $\qpropp = \big(\qpropp(m)\big)_{m\in \Z}$ as a restricted version of PROPs by putting the number of ``outputs'' to $1$. Roughly speaking, an element in $\qpropp(m)$, where $m\in \Z$, gives rise to an operation with $m$ ``inputs'' and $1$ ``output''. He used them to study iterated loop spaces \bfcit{May}, \bfcit{May2}. As to the scope of operads, they suffice to describe $\eg$ algebras, but not bialgebras \bfcite{MarklShniderStasheff}{p.\ 7}. 

Eggert considers a variant of Mac Lane's PROPs \bfcit{Eggert}. First, a preoperad is defined as a strict monoidal category $\proppo$ with $(\Z,+)$ as monoid of objects. For example, $\mapo$, given by
\[
\begin{array}{r c l}
\mapo(m,n) & := & \{f : \{1,{\dots},m\} \to \{1,{\dots},n\} \,|\, f \text{ is a map}\}
\end{array}
\]
for $m,n\in \Z$, is a preoperad. The opposite category of $\mapo$, denoted as $\mapoop$, is also a preoperad.

Then, an Eggert-operad is a tuple $\oppo$ consisting of a preoperad $\proppo$ and a morphism of preoperads $\mathfrak{p}_0 : \mapoop \to \proppo$ satisfying a Braiding Property and a Branching Property; \cf definition \ref{VER1_DEFI_OPERAD}. 

So in contrast to Mac Lans's PROPs, where to each bijective map $\{1,{\dots},m\} \to \{1,{\dots},m\}$ an element of $\rproppo(m,m)$ is attached, in an Eggert-operad to each arbitrary map $\{1,{\dots},m\} \to \{1,{\dots},n\}$ an element of $\proppo(m,n)$ is attached, where $m,n\in \Z$.

For example, given a set $M$, we can define the Eggert-operad $\Eno(M) = (\eno(M),\mathfrak{e}_0)$ with the underlying preoperad $\eno(M)$ given by
\[
\begin{array}{r c l}
\eno(M)(m,n) & := & \{f : M^{\times m} \to M^{\times n} \,|\, f \text{ is a map}\}
\end{array}
\]
for $m,n\in \Z$ and a certain morphism $\mathfrak{e}_0 : \mapoop \to \eno(M)$ of preoperads. 

Given an Eggert-operad $\oppo$, a $\oppo$-algebra is then defined as a tuple $(M,\varphi_0)$ consisting of a set $M$ and a morphism $\varphi : \oppo \to \Eno(M)$. 

A biindexed set is just a tuple of sets $X = \big(X(m,n)\big)_{m,n\in\Z}$. We construct the free Eggert-operad $\Freeox$ and show that it satisfies a universal property. We point out how to construct Eggert-operads with generators and relations using the free Eggert-operad.

To illustrate the concept, we construct the Eggert-operad $\Groupfp$ and prove that every group is a $\Groupfp$-algebra and vice versa.

We prove operadic versions of some elementary facts about groups such as the uniqueness of the inverse element.

It is a well-known fact that to define a group it suffices to require the associativity of the multiplication, left-neutrality and the existence of left-inverse elements. We provide an operadic version of this fact.

\section{Preliminaries}

\subsection{Basic notations}

\bnota \blue{}
\begin{itemize}
\item We write $\Z := \{m\in \mathbb{Z}\,|\, m\geq 0\}$ and $[m,n]  :=  \{i\in \Z \,|\, m\leq i\leq n\}$ for $m,n\in \Z$.
\item Let $m\in\Z$ and $M$ be a set. We write $M^{\times m} := \{(x_1,{\dots},x_m)\,|\, x_i\in M \text{ for } i\in [1,m]\}$. In particular we have $M^{\times 0} = \{()\}$. We identify $M^{\times 1} = M$.
\item We denote the composition of morphisms from left to right; $\ie$ given morphisms $f : X \to Y$ and $g : Y \to Z$ in a category, we will write $fg = f \cdot g : X \to Z$ for the composite of $f$ and $g$.
\item In accordance with the notation of Eggert, we add an index $0$ in the notation of preoperads, Eggert-operads etc.\ to indicate that we work over the category of sets.
\end{itemize}
\enota

\bdefi \label{VER1_DEFI_BIINDEXED_SETS} 
We define a \textit{biindexed set} as a tuple $X = (X,\so_X,\ta_X)$ consisting of a set $X$, the \textit{source map} $\so = \so_X : X \to \Z$ and the \textit{target map} $\ta = \ta_X : X \to \Z$. Further we define $X(m,n)  :=  \{x\in X\,|\, x\so_X = m,x\ta_X=n\}$ for $ m,n\in\Z$. We will often denote a biindexed set as $X = (X,\so_X,\ta_X)$ or as $\big(X(m,n)\big)_{m,n\in\Z} = (X,\so_X,\ta_X)$; \cf \bfcite{Eggert}{Def.\ 2.1}.

Let $X$, $Y$ be biindexed sets. We define $(X\btimes Y)(m,n) := X(m,n)\times Y(m,n)$ for $m,n\in \Z$. Further we define the biindexed set $X \btimes Y := \big((X\btimes Y)(m,n)\big)_{m,n\in \Z}$. So $(x,y)\in X\btimes Y$ has source $(x,y)\so := x\so = y\so$ and target $(x,y)\ta := x\ta = y\ta$; \cf \bfcite{Eggert}{Def.\ 2.5}.

Let $X$ be a biindexed set. We define $X^{\btimes 2} := X \btimes X$.

Let $X$, $Y$ be biindexed sets and $g : X \to Y$ a map. We call $g$ a \textit{morphism of biindexed sets}, or \textit{biindexed map}, if $(xg)\so_Y = x\so_X$ and $(xg)\ta_Y = x\ta_X$ for $x\in X$; \cf \bfcite{Eggert}{Def.\ 2.2}.
\edefi

\bdefi[\bfcite{Eggert}{Def.\ 2.6, Rem.\ 2.8}]
We define a \textit{preoperad}, or \textit{set-preoperad}, as a strict monoidal category $\proppo$ with object monoid $\catob_{\proppo}  =  (\Z,+)$ and $(\boxtimes)  =  (\boxtimes_{\proppo})  :  \proppo \times \proppo  \to  \proppo$ as monoidal product.
\edefi

\bbei[\bfcite{Eggert}{Def.\ 2.57}] \label{VER1_BEI_MAPO}
Let $m,n\in \Z$. We define 
\[
\begin{array}{r c l}
\mapo(m,n) & := & \{f:[1,m] \to [1,n] \,|\, f\text{ is a map}\}.
\end{array}
\]
In this way we obtain a category $\mapo$ using identities $\id_m := \id_{m,\mapo} := \id_{[1,m]} : [1,m] \to [1,m]$ for $m\in \Z$ and composition of maps.

Let $m,m',n,n'\in\Z$ and $f : [1,m] \to [1,n]$, $f' : [1,m'] \to [1,n']$ be maps. We define the map
\[
\begin{array}{r c c c c c l}
f\boxtimes f' & \stackrel{\text{short}}{=} & f\boxtimes_{\mapo} f' & : & [1,m+m'] & \to & [1,n+n']\\
 & & & & i &\mapsto & \begin{cases}
if &\text{ if } i\in[1,m]\\
n + (i-m)f'  &\text{ if }  i\in[m+1,m+m'].
\end{cases}
\end{array}
\]
This yields a functor
\[
\begin{array}{r c c c c c l}
(\boxtimes) & \stackrel{\text{short}}{=} & (\boxtimes_{\mapo}) & : & \mapo\times \mapo & \to & \mapo \\
 & & & & \big((m,m') \xrightarrow{(f,f')} (n,n')\big) & \mapsto & (m+m'  \xrightarrow{f\boxtimes f'} n+n').
\end{array}
\]
%Further we define the multiplication
%\[
%\begin{array}{r c c c l}
%(\boxtimes_{\mapo})_{m,m',n,n'} & : & \mapo(m,n) \times \mapo(m',n') & \to & \mapo(m+m',n+n') \\
% & & (f,f') & \mapsto & f\boxtimes f'.
%\end{array}
%\]
%We define the biindexed set
%\[
%\begin{array}{r c l}
%\mapo & := & \big(\mapo(m,n)\big)_{m,n\in\Z}.
%\end{array}
%\]
Endowed with the multiplication $(\boxtimes)$, the category $\mapo$ becomes a preoperad, again denoted as $\mapo$.
\ebei

\bbei[\Cf \bfcite{Eggert}{Def.\ 2.13}] \label{VER1_BEI_MAPOOP}
We consider the opposite category $\mapoop$. We write
\[
\begin{array}{r c l}
\mapoop(n,m) & = & \{f\op \,|\, f\in \mapo(m,n)\}
\end{array}
\]
for $n,m\in\Z$. So we have $(f\cdot g)\op = g\op \cdot f\op$ for $f\in \mapo(m,n)$ and $g\in\mapo(n,o)$, where $m,n,o\in\Z$. We endow $\mapoop$ with the monoidal product
\[
\begin{array}{r c c c l}
f\op \boxtimes {f'}\op & := & (f\boxtimes f')\op & \in & \mapoop(n+n',m+m')
\end{array}
\]
for $f\in \mapo(m,n)$ and $f'\in \mapo(m',n')$, where $m,n,m',n'\in \Z$. Then $\mapoop$ is a preoperad.
\ebei

\bbei[\bfcite{Eggert}{Def.\ 2.55}]
Let $m,n\in\Z$. We define 
\[
\begin{array}{r c l}
\eno(M)(m,n) & := & \{f : M^{\times m} \to M^{\times n} \,|\, f\text{ is a map}\}.
\end{array}
\]
In this way we obtain a category $\eno(M)$ using identities $\id_{m,\eno(M)} := \id_{M^{\times m}} : M^{\times m} \to M^{\times m}$ for $m\in \Z$ and composition of maps.

%Let $m,n,o\in \Z$. We define the composition
%\[
%\begin{array}{r c c c l}
%(\cdot_{\,\eno(M)})_{m,n,o} & : & \eno(M)(m,n) \times \eno(M)(n,o) & \to & \eno(M)(m,o) \\
% & & (f,g) & \mapsto & f\cdot g.
%\end{array}
%\]
Let $(x_1,{\dots},x_m)\in M^{\times m}$ and $(y_1,{\dots},y_n)\in M^{\times n}$ be tuples. We define the tuple 
\[
\begin{array}{r c c c l}
(x_1,{\dots},x_m)\times (y_1,{\dots},y_n) & := & (x_1,{\dots},x_m,y_1,{\dots},y_n) & \in & M^{\times (m+n)}
\end{array}
\]
Let $m,m',n,n'\in \Z$ and $f : M^{\times m} \to M^{\times n}$, $f' : M^{\times m'} \to M^{\times n'}$ be maps. We define the map 
\[
\begin{array}{r c c c c c l}
f\boxtimes f' & \stackrel{\text{short}}{=} & f\boxtimes_{\eno(M)} f' & : &  M^{\times (m+m')} & \to & M^{\times (n+n')} \\
 & & & & (x_1,{\dots} ,  x_{m+m'}) & \mapsto & (x_1,{\dots} , x_{m})f \times (x_{m+1}, {\dots} , x_{m+m'})f'.
\end{array}
\]
This yields a functor
\[
\begin{array}{r c c c c c l}
(\boxtimes) & \stackrel{\text{short}}{=} & (\boxtimes_{\eno(M)}) & : & \eno(M)\times \eno(M) & \to & \eno(M) \\
 & & & & \big((m,m') \xrightarrow{(f,f')} (n,n')\big) & \mapsto & (m+m'  \xrightarrow{f\boxtimes f'} n+n').
\end{array}
\]
%We define the biindexed set
%\[
%\begin{array}{r c l}
%\eno(M) & := & \big(\eno(M)(m,n)\big)_{m,n\in\Z}.
%\end{array}
%\]
Endowed with the multiplication $(\boxtimes)$, the category $\eno(M)$ becomes a preoperad, again denoted as $\eno(M)$.
\ebei

\bdefi[\bfcite{Eggert}{Def.\ 2.14}] \label{VER1_DEFI_MORPHISM_OF_PREOPERADS}
We define a \textit{morphism of preoperads} as a strict monoidal functor $\varphi_0  :  \proppo  \to  \qproppo$ between preoperads $\proppo$ and $\qproppo$ that acts identically on the objects. In particular $\varphi_0$ is a morphism of biindexed sets.
\edefi

\subsection{Eggert-operads}

\bdefi \label{VER1_DEFI_SE_UND_HE} 
Let $m,m'\in \Z$. We define the permutation
\[
\begin{array}{r c c c l}
\se_{m,m'} & : & [1,m+m'] & \to & [1,m+m'] \\
 & & i & \mapsto & \begin{cases}
 i + m' & \text{ if } 1\leq i\leq m \\
 i-m & \text{ if } m+1\leq i \leq m+m';
 \end{cases}
\end{array}
\]\cf \bfcite{Eggert}{Def.\ 6.1}. Further we define the map
\[
\begin{array}{r c c c l}
\he_{m,m'} & : &  [1,mm'] & \to & [1,m'] \\
 & & i &\mapsto & \underline{i},
\end{array}
\]
where $\underline{i}\in[1,m']$ is the unique integer such that $m'$ divides $i-\underline{i}$; \cf \bfcite{Eggert}{Def.\ 6.2}.

Let $a\in \Z$ and $\proppo$ be a preoperad. Let $x\in \proppo$. We define recursively $x^{\boxtimes 0} := \id_{0,\proppo}$ and 

$x^{\boxtimes a} := x^{\boxtimes (a-1)} \boxtimes x$ if $a\geq 1$; \cf \bfcite{Eggert}{Def.\ 2.6 (mc2)}.
\edefi

\bdefi[\bfcite{Eggert}{Def.\ 6.3}] \label{VER1_DEFI_OPERAD}
We define an \textit{(Eggert-)operad} as a tuple $(\proppo,\mathfrak{p}_0)$ consisting of a preoperad $\proppo$ and a morphism of preoperads $\mathfrak{p}_0  :  \mapoop  \to  \proppo $, called \textit{structure morphism}, such that the following axioms hold:

\underline{Braiding Property}: For $x,x'\in \proppo$, we have 
\[
\begin{array}{r c l}
\se_{x\so,x'\so}\op\mathfrak{p}_0 \cdot (x \boxtimes x') & = & (x' \boxtimes x) \cdot \se_{x\ta,x'\ta}\op\mathfrak{p}_0.
\end{array}
\]
\underline{Branching Property}: For $a\in \Z$ and $x\in \proppo$, we have 
\[
\begin{array}{r c l}
\he_{a,x\so}\op\mathfrak{p}_0 \cdot x^{\boxtimes a} & = & x \cdot \he_{a,x\ta}\op\mathfrak{p}_0.
\end{array}
\]
\edefi

\bbem \label{VER1_BEM_VERGLEICH_PROP_UND_EGGERT_OPERADE}
Definition \ref{VER1_DEFI_OPERAD} is a variant of the notion of PROPs, \cf \bfcite{MacLane}{p.\ 97}. In comparison, the data of a PROP include the images of the bijective maps under $\mathfrak{p}_0$. The Braiding Property corresponds to equation (4) in \bfcite{MacLane}{p.\ 97}.
\ebem

\bbei \label{VER1_BEI_OPERADE_ENO(M)}  
Let $M$ be a set. Let $m,n\in \Z$ and $f\in \mapo(m,n)$. We define the map
\[
\begin{array}{r c c c l}
f\op\mathfrak{e}_0 & : & M^{\times n} & \to & M^{\times m} \\
 & & (x_1,{\dots},x_n) & \mapsto & (x_{1f},{\dots},x_{mf});
\end{array}
\] 
\cf \bfcite{Eggert}{Def.\ 2.61}. Further we define the morphism of preoperads
\[
\begin{array}{r c c c l}
\mathfrak{e}_0 & : & \mapoop & \to & \eno(M) \\
 & & f\op & \mapsto & f\op\mathfrak{e}_0;
\end{array}
\]
\cf \bfcite{Eggert}{Lem.\ 2.62}.
Then $\mathfrak{e}_0$ satisfies the Braiding Property and Branching Property, so
\[
\begin{array}{r c l}
\Eno(M) & := & (\eno(M),\mathfrak{e}_0)
\end{array}
\]
is an operad; \cf \bfcite{Eggert}{Ex.\ 6.5}.
\ebei

\bdefi[\bfcite{Eggert}{Def.\ 6.13}]
Let $\oppo$ and $\qoppo$ be operads. We define a \textit{morphism of operads} $\varphi_0 : \oppo \to \qoppo$ as a morphism of preoperads $\varphi_0 : \proppo \to \qproppo$ such that $\mathfrak{p}_0 \cdot \varphi_0 = \mathfrak{q}_0$.
\edefi

\bdefi[\bfcite{Eggert}{Def.\ 6.21}] \label{VER1_DEFI_ALGEBRA}
Let $\oppo$ be an operad. We define a \textit{$\oppo$-algebra} as a tuple $(M,\varphi_0)$ consisting of a set $M$ and a morphism of operads $\varphi_0 : \oppo \to \Eno(M)$.
\edefi

\section{The free Eggert-operad $\Freeox$}

Throughout this section, let $X = (X,\so_X,\ta_X)$ be a biindexed set.

\subsection{Words over the biindexed set $X$}

\bdefi \label{VER1_DEFI_WORDS}
We define a \textit{word} (over $X$) as a $(2k+1)$-tuple $w =  (a_1,{\dots},a_{2k+1})$, where 
\[
\begin{array}{r c l l }
\vsp k & \in & \Z\\
\vsp a_{2i+1} & =: & f_i\op \in \mapoop  & \text{ for } i\in [0,k] \\
\vsp a_{2i} & =: & \big(l_i,x_i,r_i\big)\in \{\big(l,x,r\big)\,| \,l,r\in\Z, x\in X\} &\text{ for } i\in [1,k]
\end{array}
\]
such that 
\[
\begin{array}{r c c c l l }
\vsp  &  & f_i\op\ta_{\mapoop} & = & l_{i+1} + x_{i+1}\so_{X} + r_{i+1} & \text{ for } i\in [0,k-1]\\
\vsp l_i + x_i\ta_{X} + r_i & = & f_i\op\so_{\mapoop} & &  & \text{ for } i\in [1,k].
\end{array}
\]
Further we call $k$ the \textit{length} of the word $w$ and denote it as $|w|  :=  k  \in  \Z$. We denote the set of words over $X$ as $\Wordso$.

Let $w =  (a_1,{\dots},a_{2k+1})\in \Wordso$ be a word. We define its source $w\so := a_1\so$ and its target $w\ta := a_{2k+1}\ta$.

Endowed with the source and target defined above, $\Wordso$ becomes a biindexed set, again denoted as $\Wordso$.

We endow $\Wordso$ with the composition
\[
\begin{array}{r c l}
(a_1,{\dots},a_{2k+1}) \cdot (a'_1,a'_2,{\dots},a'_{2k'+1}) & := & (a_1,{\dots},a_{2k}, a_{2k+1} \cdot a'_1 ,a'_2,{\dots},a'_{2k'+1}),
\end{array}
\]
where $(a_1,{\dots},a_{2k+1}), (a'_1,a'_2,{\dots},a'_{2k'+1}) \in \Wordso$ are words such that $(a_1,{\dots},a_{2k+1})\ta = (a'_1,a'_2,{\dots},a'_{2k'+1})\so$. Note that $a_{2k+1} \cdot a'_1$ is the composite formed in $\mapoop$. 

With the identities $\id_{m,\Wordso} = (\id_m\op) \in \Wordso(m,m)$ for $m\in \Z$ and the composition defined above, the biindexed set $\Wordso$ becomes a category, again denoted as $\Wordso$.

Note that $|w \cdot w'| = |w| + |w'|$ for words $w,w'\in \Wordso$ such that $w\ta = w'\so$.

We will often denote a word $\big(f\op\big) \in \Wordso$ of length 0 as 
\[
\begin{array}{r c l}
f\op  & \stackrel{\text{short}}{=} & (f\op).
\end{array}
\]
Let $w = (a_1,{\dots},a_{2k+1})\in\Wordso$ be a word. Note that there exist unique $f_i\in \mapo$ such that $a_{2i+1} = f_i\op$ for $i\in[0,k]$. Note that there exist unique $l_i,r_i\in \Z$ and $x_i\in X$ such that $a_{2i} = \big(l_i,x_i,r_i\big)$ for $i\in[1,k]$. We define the words
\[
\begin{array}{r c c c l}
w_i & := & (\id_{l_i+x_i\so_X+r_i}\op,\big(l_i,x_i,r_i\big),\id_{l_i+x_i\ta_X+r_i}\op) & \in & \Wordso
\end{array}
\]
of length 1 for $i\in[1,k]$. Further we define the word
\[
\begin{array}{r c l}
f_0\op \cdot w_1 \cdot f_1\op \cdot {\dots} \cdot w_j \cdot f_j\op & := & \begin{cases}
f_0\op & \text{ if } j=0 \\
(f_0\op \cdot w_1 \cdot f_1\op \cdot {\dots} \cdot w_{j-1} \cdot f_{j-1}\op) \cdot w_j \cdot f_j\op & \text{ if } 1 \leq j\leq k
\end{cases}
\end{array}
\]
for $j\in[0,k]$. Note that
\[
\begin{array}{r c c c l}
f_0\op \cdot w_1 \cdot f_1\op \cdot {\dots} \cdot w_k \cdot f_k\op & = & w & \in & \Wordso.
\end{array}
\]
We call $f_0\op \cdot w_1 \cdot f_1\op \cdot {\dots} \cdot w_k \cdot f_k\op $ the \textit{standard decomposition} of $w$. 
%Note that given two standard decompositions of a word
%\[
%\begin{array}{r c l}
%f_0\op \cdot w_1 \cdot f_1\op \cdot {\dots} \cdot w_k \cdot f_k\op & = & {f_0'}\op \cdot w_1' \cdot {f_1'}\op \cdot {\dots} \cdot w'_{k'} \cdot {f'_{k'}}\op,
%\end{array}
%\]
%we have $k = k'$, $f_i\op = {f_i'}\op$ for $i\in [0,k]$ and $w_i = w_i'$ for $i\in[1,k]$. So the standard decomposition of a word is unique.

Let $q,p,l,r\in \Z$, $x\in X$ and $f \in \mapo$. We define
\[
\begin{array}{r c l}
f\op\triangleright p & := & f\op \boxtimes \id_p\op\\
q\triangleleft f\op & := & \id_q\op \boxtimes \, f\op\\
q\triangleleft f\op \triangleright p & := & \id_q\op \boxtimes \, f\op \boxtimes \id_p\op\\
\big(l,x,r\big)\triangleright p & := & \big(l,x,r+p\big)\\
q \triangleleft \big(l,x,r\big) & := & \big(q+l,x,r\big) \\
q \triangleleft \big(l,x,r\big)\triangleright p & := & \big(q+l,x,r+p\big).
\end{array}
\]
Let $q,p,m,n\in \Z$ and $w = \big(a_1,{\dots},a_{2k+1}\big)\in\Wordso(m,n)$ be a word. We define
\[
\begin{array}{r c l c l}
w \triangleright p & := & \big(a_1\triangleright p,{\dots},a_{2k+1}\triangleright p\big) & \in & \Wordso(m+p,n+p)\\
q\triangleleft w & := & \big(q \triangleleft a_1,{\dots},q \triangleleft a_{2k+1}\big) & \in & \Wordso(q+m,q+n)\\
q\triangleleft w \triangleright p & := & \big(q\triangleleft a_1\triangleright p,{\dots},q\triangleleft a_{2k+1}\triangleright p\big) & \in & \Wordso(q+m+p,q+n+p).
\end{array}
\]
We will often denote a word of the form $w = (\id_{l+x\so_X+r}\op,\big(l,x,r\big),\id_{l+x\ta_X+r}\op) \in \Wordso$ as
\[
\begin{array}{r c l}
\big(l,x,r\big) & \stackrel{\text{short}}{=} & (\id_{l+x\so_X+r}\op,\big(l,x,r\big),\id_{l+x\ta_X+r}\op).
\end{array}
\] 
\edefi

\bbem \label{VER1_BEM_FORMELN_FÜR_TRIANGLE}
Let $p,q\in \Z$ and $w,w'\in \Wordso$ be words such that $w\ta = w'\so$. The following equations hold.
\[
\begin{array}{r c l c r c l}
0 \triangleleft w & = & w &  & w \triangleright 0 & = & w \\
|q\triangleleft w| & = & |w| & & |w\triangleright p| & = & |w| \\
q\triangleleft (p\triangleleft w) & = & (q+p)\triangleleft w & & (w\triangleright p)\triangleright q & = & w\triangleright(p+q) \\
(q \triangleleft w)\cdot (q \triangleleft w') & = & q \triangleleft (w\cdot w') & & (w\triangleright p)\cdot (w'\triangleright p) & = & (w\cdot w')\triangleright p \\
(q \triangleleft w) \triangleright p & = & q \triangleleft w\triangleright p & & q \triangleleft (w\triangleright p) & = & q \triangleleft w\triangleright p \\
(q \triangleleft v \triangleright p)\so & = & q+v\so+p &  & (q \triangleleft v \triangleright p)\ta & = &  q+v\ta+p
\end{array}
\]
\ebem

\bdefi \label{VER1_DEFI_POTENZ_IN_WORDS}
Let $a,m,n\in \Z$ and $w\in \Wordso(m,n)$ be a word. We define recursively
\[
\begin{array}{r c c c l}
w^{\boxtimes a} & := & \left. \begin{cases}
\id_0 & \text{ if } a=0 \\
(w^{\boxtimes (a-1)} \triangleright w\so) \cdot (w^{\boxtimes (a-1)}\ta \triangleleft w) & \text{ if } a\geq 1
\end{cases} \right\} & \in & \Wordso(am,an).
\end{array}
\]
This is a provisional notation. Recall that we have not introduced a monoidal product on the category $\Wordso$. The connection between $w^{\boxtimes a}$ and an actual tensor power will be clarified in definition \ref{VER1_DEFI_FREEOPREX} and lemma \ref{VER1_LEM_VORBEREITUNG_STRUKTURMORPHISMUS_VON_FREEOX}.1 below.
\edefi

\bdefi \label{VER1_DEFI_ELEMENTARE_ÄQUIVALENZEN}
We define the following subsets of $\big(\Wordso\big)^{\btimes 2}$; \cf definition \ref{VER1_DEFI_BIINDEXED_SETS}. 
\[
\begin{array}{r l l}
M^1_{0,X}  :=  \{ (w,w') \in \big(\Wordso\big)^{\btimes 2} \,| &\hspace{-2ex} w = (v\triangleright v'\so) \cdot (v\ta \triangleleft v'),\,  w' = (v\so \triangleleft v') \cdot (v\triangleright v'\ta) & \hspace{-2ex} , \\
 & \hspace{-2ex} v,v'\in \Wordso,|v|,|v'| \leq 1 & \hspace{-2ex} \} \\
M^2_{0,X}  :=  \{ (w,w') \in  \big(\Wordso\big)^{\btimes 2} \,| & \hspace{-2ex}  w = q\triangleleft \big(\se_{v\so,a}\op \cdot (v\triangleright a)\big) \triangleright p,\, w' = q\triangleleft \big((a\triangleleft v) \cdot \se_{v\ta,a}\op\big) \triangleright p & \hspace{-2ex}, \\
 & \hspace{-2ex} a,q,p \in \Z, v\in \Wordso, |v|\leq 1 & \hspace{-2ex} \} \\
M^3_{0,X}  :=  \{ (w,w') \in  \big(\Wordso\big)^{\btimes 2} \,| & \hspace{-2ex}  w = q\triangleleft \big(\se_{a,v\so}\op \cdot (a\triangleleft v)\big) \triangleright p,\,  w' = q\triangleleft \big((v\triangleright a) \cdot \se_{a,v\ta}\op\big) \triangleright p &  \hspace{-2ex} , \\
 & \hspace{-2ex}  a,q,p \in \Z, v\in \Wordso, |v|\leq 1 & \hspace{-2ex} \} \\
M^4_{0,X}  :=  \{ (w,w') \in  \big(\Wordso\big)^{\btimes 2}  \,| & \hspace{-2ex} w = q\triangleleft (\he_{a,v\so}\op \cdot v^{\boxtimes a}) \triangleright p,\, w' = q\triangleleft (v \cdot \he_{a,v\ta}\op) \triangleright p & \hspace{-2ex} , \\
 & \hspace{-2ex}  a,q,p \in \Z, v\in \Wordso, |v|\leq 1 & \hspace{-2ex} \}
\end{array}
\]
Further we define 
\[
\begin{array}{r c c c l}
M_{0,X} & := & \displaystyle \bigcup_{i=1}^4 M^i_{0,X} & \subseteq & \big(\Wordso\big)^{\btimes 2}.
\end{array}
\]
We write
\[
\begin{array}{r c l}
w\stackrel{i.}{\sim} w' & :\Longleftrightarrow & (w,w') \in M^i_{0,X}
\end{array}
\]
for $i\in [1,4]$ and
\[
\begin{array}{r c l}
w\sim w' & :\Longleftrightarrow & (w,w') \in M_{0,X}.
\end{array}
\]
Note that $w\so= w'\so$ and $w\ta = w'\ta$ for $(w,w')\in M_{0,X}$.

We define $N_{0,X} \subseteq \big(\Wordso\big)^{\btimes 2}$ as the congruence relation generated by $M_{0,X}$. Note that $N_{0,X}$ is the equivalence relation on $\Wordso$ generated by 
\[
\begin{array}{r l}
M_{0,X}'  :=  \{(uav,ubv) \in \big(\Wordso\big)^{\btimes 2} \,| & \hspace{-2ex} \big((a,b)\in M_{0,X} \vee (b,a)\in M_{0,X}\big), \\
 & \hspace{-2ex} u,v\in \Wordso, u\ta = a\so=b\so, a\ta=b\ta=v\so\}.
\end{array}
\]
We write
\[
\begin{array}{r c l}
w \approx w' & :\Longleftrightarrow & (w,w') \in N_{0,X}.
\end{array}
\]
Note that $w\so= w'\so$ and $w\ta = w'\ta$ for $(w,w')\in N_{0,X}$.
\edefi

\blem \label{VER1_LEM_VORBEREITUNG_AUF_DEFINITION_FREEOPREX} \blue{} 
\begin{enumerate}
\item \textit{Let $d,b\in \Z$ and $w,w'\in \Wordso$ be words such that $w\sim w'$. Then we have
\[
\begin{array}{r c l}
d \triangleleft w \triangleright b & \sim & d \triangleleft w' \triangleright b.
\end{array}
\]}
\item \textit{Let $d,b\in \Z$ and $w,w'\in \Wordso$ be words such that $w\approx w'$. Then we have
\[
\begin{array}{r c l}
d \triangleleft w \triangleright b & \approx & d \triangleleft w' \triangleright b.
\end{array}
\]}
\item \textit{Let $u,u',v,v'\in \Wordso$ be words such that $u\approx u'$ and $v\approx v'$. Then we have
\[
\begin{array}{r c l}
(u\triangleright v\so) \cdot (u\ta \triangleleft v) & \approx & (u'\triangleright v'\so) \cdot (u'\ta \triangleleft v').
\end{array}
\]}
\end{enumerate}
\bpro Ad 1: In all four cases $w\stackrel{i.}{\sim} w'$ of definition \ref{VER1_DEFI_ELEMENTARE_ÄQUIVALENZEN}, where $i\in [1,4]$, this  follows from remark \ref{VER1_BEM_FORMELN_FÜR_TRIANGLE}.

Ad 2: We may choose words $z_0, {\dots}, z_k \in \Wordso$ such that $w=z_0$, $z_k = w'$ and $(z_i,z_{i+1})\in M_{0,X}'$ for $i\in [0,k-1]$. By 1. and remark \ref{VER1_BEM_FORMELN_FÜR_TRIANGLE} we conclude that $(d\triangleleft z_i \triangleright b, d\triangleleft z_{i+1} \triangleright b) \in M_{0,X}'$ for $i\in [0,k-1]$, hence $d \triangleleft w \triangleright b  \approx  d \triangleleft w' \triangleright b$.

Ad 3: By 2. and $v\so = v'\so$ we have $(u \triangleright v\so) \approx (u' \triangleright v'\so)$. By 2. and $u\ta = u'\ta$ we have $(u\ta \triangleleft v) \approx (u'\ta \triangleleft v')$. Since $N_{0,X}$ is a congruence relation, we conclude that $(u\triangleright v\so) \cdot (u\ta \triangleleft v) \approx (u'\triangleright v'\so) \cdot (u'\ta \triangleleft v')$.
\epro
\elem

\blem \label{VER1_LEM_FREEOPREX_INTERCHANGING_RULE}
\textit{Let $w,w'\in \Wordso$ be words. Then we have
\[
\begin{array}{r c l}
(w \triangleright w'\so) \cdot (w\ta \triangleleft w') & \approx & (w\so \triangleleft w') \cdot (w \triangleright w'\ta).
\end{array}
\]
}
\bpro We consider particular cases:

\textit{Case 1: We have $|w|, |w'| \leq 1$.}

Then we have $(w \triangleright w'\so) \cdot (w\ta \triangleleft w') \stackrel{1.}{\sim} (w\so \triangleleft w') \cdot (w \triangleright w'\ta)$, hence $(w \triangleright w'\so) \cdot (w\ta \triangleleft w')  \approx  (w\so \triangleleft w') \cdot (w \triangleright w'\ta)$; \cf definition \ref{VER1_DEFI_ELEMENTARE_ÄQUIVALENZEN}.

\textit{Case 2: We have $|w'| \leq 1$.}

Let $k := |w|$. If $k\leq 1$, Case 1 finishes the proof. Otherwise let $w = f_0\op \cdot w_1 \cdot f_1\op \cdot {\dots} \cdot w_k \cdot f_k\op$ be the standard decomposition. We define the word $u := f_0\op \cdot w_1 \cdot  f_1\op \cdot {\dots} \cdot w_{k-1} \cdot f_{k-1}\op \in \Wordso$ of length $k-1$ and the word $v := w_k \cdot f_k\op \in \Wordso$ of length $1$. Note that $w = u \cdot v$, $w\so = u\so$, $u\ta = v\so$, $v\ta = w\ta$. Using induction over $k$, we conclude that
\[
\begin{array}{c l c l}
 & (w \triangleright w'\so) \cdot (w\ta \triangleleft w') & = & \big((u\cdot v) \triangleright w'\so\big) \cdot (v\ta \triangleleft w')\\
 \stackrel{\ref{VER1_BEM_FORMELN_FÜR_TRIANGLE}}{=} & (u \triangleright w'\so) \cdot (v \triangleright w'\so) \cdot (v\ta \triangleleft w') & \stackrel{\text{Case 1}}{\approx} & (u \triangleright w'\so) \cdot (v\so \triangleleft w') \cdot (v \triangleright w'\ta) \\
 = & (u \triangleright w'\so) \cdot (u\ta \triangleleft w') \cdot (v \triangleright w'\ta)  & \stackrel{\IV}{\approx} & (u\so \triangleleft w') \cdot (u \triangleright w'\ta) \cdot (v \triangleright w'\ta) \\
 \stackrel{\ref{VER1_BEM_FORMELN_FÜR_TRIANGLE}}{=} & (u\so \triangleleft w') \cdot \big((u\cdot v) \triangleright w'\ta\big)  & = & (w\so \triangleleft w') \cdot (w \triangleright w'\ta).
\end{array}
\]

\textit{Case 3: We have $|w'| \geq 2$.}

Let $k' := |w'|$. Let $w' = {f_0'}\op \cdot w_1' \cdot {f_1'}\op \cdot {\dots} \cdot w_{k'}' \cdot {f_{k'}'}\op$ be the standard decomposition. We define the word $u' := {f_0'}\op \cdot w'_1 \cdot  {f_1'}\op \cdot {\dots} \cdot w'_{k'-1} \cdot {f'_{k'-1}}\hspace{-3ex}\op \in \Wordso$ of length $k'-1$ and the word $v' := w'_{k'} \cdot {f'_{k'}}\op \in \Wordso$ of length $1$. Note that $w' = u' \cdot v'$, $w'\so = u'\so$, $u'\ta = v'\so$, $v'\ta = w'\ta$. Using induction over $k'$, we conclude that
\[
\begin{array}{c l c l}
 & (w \triangleright w'\so) \cdot (w\ta \triangleleft w') & = & (w \triangleright u'\so) \cdot \big(w\ta \triangleleft (u' \cdot v')\big)\\
 \stackrel{\ref{VER1_BEM_FORMELN_FÜR_TRIANGLE}}{=} & (w \triangleright u'\so) \cdot (w\ta \triangleleft u') \cdot (w\ta \triangleleft v') & \stackrel{\IV}{\approx} & (w\so \triangleleft u') \cdot (w \triangleright u'\ta) \cdot (w\ta \triangleleft v') \\
 = & (w\so \triangleleft u') \cdot (w \triangleright v'\so) \cdot (w\ta \triangleleft v')  & \stackrel{\text{Case 2}}{\approx} & (w\so \triangleleft u') \cdot (w\so \triangleleft v') \cdot (w \triangleright v'\ta) \\
 \stackrel{\ref{VER1_BEM_FORMELN_FÜR_TRIANGLE}}{=} & \big(w\so \triangleleft (u' \cdot v')\big) \cdot (w \triangleright v'\ta)  & = & (w\so \triangleleft w') \cdot (w \triangleright w'\ta).
\end{array}
\]
\epro
\elem

\subsection{Construction of $\Freeox$}

\bdefi \label{VER1_DEFI_FREEOPREX}
We define the factor category of $\Wordso$ modulo $N_{0,X}$
\[
\begin{array}{r c l}
\Freeoprex & := & \Wordso / N_{0,X}.
\end{array}
\]
This construction is also known as quotient category. 

We shall now endow the category $\Freeoprex$ with a multiplication in order to obtain a preoperad, again denoted as $\Freeoprex$.
Let $m,n,m',n'\in \Z$ and $w\in\Wordso(m,n)$, $w'\in \Wordso(m',n')$ be words. We define
\[
\begin{array}{r c c c l}
[w] \boxtimes [w'] & := & [(w \triangleright w'\so) \cdot (w\ta \triangleleft w')] & \in & \Freeoprex(m+m',n+n'),
\end{array}
\]
which is well-defined by lemma \ref{VER1_LEM_VORBEREITUNG_AUF_DEFINITION_FREEOPREX}.3. We shall verify in lemma \ref{VER1_LEM_FREEOPREX_IST_PRÄOPERADE} that $\Freeoprex$ is in fact a preoperad. The upper index ``pre'' indicates that $\Freeoprex$ will be the underlying preoperad of an operad $\Freeox$, which we will introduce in definition \ref{VER1_DEFI_FREEOX}.
\edefi

\bbem \label{VER1_BEM_RECHENREGELN_FÜR_WÖRTER_DER_LÄNGE_NULL}
Let $m,n,o\in \Z$ and $f\in \mapo(m,n)$ and $g\in \mapo(n,o)$. Note that in $\Freeoprex$ we have $[g\op] \cdot [f\op] = [(f\cdot g)\op]$.

Let $m,n,m',n'\in \Z$ and $f\in \mapo(m,n)$ and $f'\in \mapo(m',n')$. Note that in $\Freeoprex$ we have $[f\op] \boxtimes [{f'}\op] = [(f\op \triangleright {f'}\op\so)\cdot (f\op\ta \triangleleft {f'}\op)] = [(f\op \boxtimes \id\op_{{f'}\op\so}) \cdot (\id\op_{f\op\ta}\boxtimes {f'}\op)] = [f\op\boxtimes {f'}\op]$.
\ebem

\blem \label{VER1_LEM_FREEOPREX_IST_PRÄOPERADE}
\textit{The category $\Freeoprex$ together with the multiplication $(\boxtimes)$ given by definition \ref{VER1_DEFI_FREEOPREX} is a preoperad.}
\bpro
In the following proof we will make repeated use of remark \ref{VER1_BEM_FORMELN_FÜR_TRIANGLE}.

Functoriality of $(\boxtimes)$ on identities: Let $m,m'\in \Z$. Then in $\Freeoprex$ we have 
\[
\begin{array}{c l c l}
 & \vsp \id_{m,\Freeoprex} \boxtimes \id_{m',\Freeoprex} & = & [\id_{m,\Wordso}] \boxtimes [\id_{m',\Wordso}] \\
 =  \vsp & [\id_{[1,m]}\op] \boxtimes [\id_{[1,m']}\op]  & \stackrel{\ref{VER1_BEM_RECHENREGELN_FÜR_WÖRTER_DER_LÄNGE_NULL}}{=} & [\id_{[1,m]}\op   \boxtimes \id_{[1,m']}\op] \\
 =  \vsp & [\id_{[1,m+m']}\op] & = & [\id_{m+m',\Wordso}] \\
 = & \id_{m+m',\Freeoprex}.
\end{array}
\]

Functoriality of $(\boxtimes)$ on composition: Let $u,u',v,v'\in \Wordso$ be words such that $u\ta = v\so$ and $u'\ta = v'\so$. Then in $\Freeoprex$ we have
\[
\begin{array}{c l c l}
 & ([u]\boxtimes [u']) \cdot ([v] \boxtimes [v']) & = & [(u\triangleright u'\so) \cdot (u\ta \triangleleft u')] \cdot [(v\triangleright v'\so) \cdot (v\ta \triangleleft v')] \\
 = & [(u\triangleright u'\so)] \cdot [(u\ta \triangleleft u') \cdot (v\triangleright v'\so)] \cdot [(v\ta \triangleleft v')] & = & [(u\triangleright u'\so)] \cdot [(v\so \triangleleft u') \cdot (v\triangleright u'\ta)] \cdot [(v\ta \triangleleft v')] \\
 \stackrel{\ref{VER1_LEM_FREEOPREX_INTERCHANGING_RULE}}{=} & [(u\triangleright u'\so)] \cdot [(v\triangleright u'\so) \cdot (v\ta \triangleleft u')] \cdot [(v\ta \triangleleft v')] & = & [(u\triangleright u'\so) \cdot (v\triangleright u'\so) \cdot (v\ta \triangleleft u') \cdot (v\ta \triangleleft v')] \\
 = & [\big((u\cdot v)\triangleright u'\so\big) \cdot \big(v\ta \triangleleft (u'\cdot v') \big)]  & = & [\big((u\cdot v)\triangleright (u' \cdot v')\so\big) \cdot \big((u \cdot v)\ta \triangleleft (u'\cdot v') \big)] \\
 = & [u\cdot v] \boxtimes [u'\cdot v'].
\end{array}
\]

Associativity of $(\boxtimes)$: Let $w,w',w''\in \Wordso$ be words. Then in $\Freeoprex$ we have
\[
\begin{array}{r c l}
([w]\boxtimes [w'])\boxtimes [w''] & = & [(w\triangleright w'\so)\cdot (w\ta \triangleleft w')] \boxtimes[w''] \\
 & = & [\Big(\big((w\triangleright w'\so)\cdot (w\ta \triangleleft w')\big) \triangleright w''\so\Big) \cdot \Big(\big((w\triangleright w'\so)\cdot (w\ta \triangleleft w')\big)\ta \triangleleft w''\Big)] \\
 & = & [\big((w \triangleright w'\so)\triangleright w''\so\big) \cdot (w\ta \triangleleft w' \triangleright w''\so) \cdot \big((w\ta \triangleleft w')\ta \triangleleft w''\big)] \\
 & = & [\big((w \triangleright (w'\so + w''\so)\big) \cdot (w\ta \triangleleft w' \triangleright w''\so) \cdot \big((w\ta + w'\ta) \triangleleft w''\big)] \\
 & = & [\big(w \triangleright (w' \triangleright w''\so)\so\big) \cdot (w\ta \triangleleft w' \triangleright w''\so) \cdot \big( w\ta \triangleleft (w'\ta \triangleleft w'')\big)] \\
 & = & [\Big(w \triangleright \big((w' \triangleright w''\so) \cdot (w'\ta \triangleleft w'')\big)\so\Big) \cdot \Big(w\ta \triangleleft \big((w' \triangleright w''\so) \cdot (w'\ta \triangleleft w'')\big)\Big)] \\
 & = & [w] \boxtimes [(w' \triangleright w''\so) \cdot (w'\ta \triangleleft w'')] \\
 & = & [w] \boxtimes ([w']\boxtimes [w'']).
\end{array}
\]

Left and right identity for $(\boxtimes)$: Let $w\in \Wordso$ be a word. Then in $\Freeoprex$ we have 
\[
\begin{array}{c l c c l}
 & \id_{0,\Freeoprex} \boxtimes\, [w] & \quad \quad  \text{and} \quad \quad \quad  & & [w] \boxtimes\id_{0,\Freeoprex}  \\
 = & [\id_{0,\Wordso}] \boxtimes [w] & & = & [w] \boxtimes [\id_{0,\Wordso}] \\
 = & [\id\op_{[1,0]}] \boxtimes [w] & & = & [w] \boxtimes [\id\op_{[1,0]}] \\
 = & [(\id_{[1,0]}\op \triangleright\, w\so) \cdot (\id_{[1,0]}\op\ta \triangleleft w)] & & = & [(w \triangleright \id_{[1,0]}\op\so) \cdot (w\ta \triangleleft \id_{[1,0]}\op)] \\
 = & [(\id_{[1,0]}\op \boxtimes \id_{[1,w\so]}\op) \cdot (0 \triangleleft w)] & & = & [(w \triangleright 0) \cdot (\id_{[1,w\ta]}\op \boxtimes \id_{[1,0]}\op)] \\
 = & [\id_{[1,w\so]}\op \cdot\, w] & & = & [w \cdot \id_{[1,w\ta]}\op]\\
 = & [w] &  & = & [w].
\end{array}
\]
\epro
\elem

\bbem \label{VER1_BEM_OXO_MIT_IDENTITÄTEN}
Let $w,w',w''\in \Wordso$ be words. Then we have 
\[
\begin{array}{r c l}
[w]\boxtimes [w']\boxtimes [w''] & = &  [\big(w\triangleright (w'\so+w''\so)\big) \cdot (w\ta \triangleleft w' \triangleright w''\so) \cdot \big((w\ta + w'\ta) \triangleleft  w''\big)].
\end{array}
\]
In particular, given $l,r\in \Z$ and $x\in X$, we have
\[
\begin{array}{r c c c l}
[\id_l\op]\boxtimes [\big(0,x,0\big)] \boxtimes [\id_r\op] & = & [\big(\id_l\op\triangleright (x\so+r)\big) \cdot (l \triangleleft \big(0,x,0\big) \triangleright r) \cdot \big((l + x\ta) \triangleleft  \id_r\op\big)] & = & [\big(l,x,r\big)].
\end{array}
\]
\ebem

\blem \label{VER1_LEM_VORBEREITUNG_SE_STRUKTURMORPHISMUS_VON_FREEOX}
\textit{Let $p\in \Z$ and $w\in \Wordso$ be a word. Then we have}
\[
\begin{array}{r c l}
[\se_{w\so,p}\op \cdot (w \triangleright p)] & = & [(p \triangleleft w) \cdot \se_{w\ta,p}\op]
\end{array}
\]
\textit{and}
\[
\begin{array}{r c l}
[\se_{p,w\so}\op \cdot (p\triangleleft w)] & = &  [(w \triangleright p) \cdot \se_{p,w\ta}\op].
\end{array}
\]
\bpro
If $k := |w| \leq 1$, this follows from $\se_{w\so,p}\op \cdot (w \triangleright p) \stackrel{2.}{\sim} (p \triangleleft w) \cdot \se_{w\ta,p}\op$ and $\se_{p,w\so}\op \cdot (p\triangleleft w) \stackrel{3.}{\sim} (w \triangleright p) \cdot \se_{p,w\ta}\op$; \cf definition \ref{VER1_DEFI_ELEMENTARE_ÄQUIVALENZEN}.

Otherwise let $w=f_0\op \cdot w_1 \cdot f_1\op \cdot {\dots} \cdot w_k \cdot f_k\op $ be the standard decomposition. We define the word $u := f_0\op \cdot w_1 \cdot  f_1\op \cdot {\dots} \cdot w_{k-1} \cdot f_{k-1}\op \in \Wordso$ of length $k-1$ and the word $v := w_k \cdot f_k\op \in \Wordso$ of length $1$. Note that $w = u \cdot v$, $w\so = u\so$, $u\ta = v\so$, $v\ta = w\ta$. Using induction over $k$, we conclude that
\[
\begin{array}{r c c c c c l}
[\se_{w\so,p}\op \cdot (w \triangleright p)] & \stackrel{\text{\ref{VER1_BEM_FORMELN_FÜR_TRIANGLE}}}{=} & [\se_{u\so,p}\op \cdot (u \triangleright p) \cdot (v \triangleright p)] & \stackrel{\IV}{=} & [(p \triangleleft u) \cdot \se_{u\ta,p}\op \cdot (v \triangleright p)] & &  \\
 & = & [(p \triangleleft u) \cdot \se_{v\so,p}\op \cdot (v \triangleright p)] 
 & \stackrel{|v| = 1}{=} & [(p \triangleleft u) \cdot (p \triangleleft v) \cdot \se_{v\ta,p}\op] & \stackrel{\text{\ref{VER1_BEM_FORMELN_FÜR_TRIANGLE}}}{=} & [(p \triangleleft w) \cdot \se_{w\ta,p}\op]
\end{array}
\]
and
\[
\begin{array}{r c c c c c l}
[\se_{p,w\so}\op \cdot (p\triangleleft w)] & \stackrel{\text{\ref{VER1_BEM_FORMELN_FÜR_TRIANGLE}}}{=} & [\se_{p,u\so}\op \cdot (p\triangleleft u) \cdot (p\triangleleft v)] & \stackrel{\IV}{=} & [(u \triangleright p) \cdot \se_{p,u\ta}\op \cdot (p\triangleleft v)] & &  \\
 & = & [(u \triangleright p) \cdot \se_{p,v\so}\op \cdot (p\triangleleft v)]
 & \stackrel{|v| = 1}{=} & [(u \triangleright p) \cdot(v \triangleright p) \cdot \se_{p,v\ta}\op] & \stackrel{\text{\ref{VER1_BEM_FORMELN_FÜR_TRIANGLE}}}{=} & [(w \triangleright p) \cdot \se_{p,w\ta}\op].
\end{array}
\]
\epro
\elem

\blem \label{VER1_LEM_VORBEREITUNG_STRUKTURMORPHISMUS_VON_FREEOX} \blue{}
\begin{enumerate}
\item \textit{Let $a\in \Z$ and $w\in \Wordso$ be a word. Then in $\Freeoprex$ we have
\[
\begin{array}{r c l}
[w]^{\boxtimes a} & = & [w^{\boxtimes a}].
\end{array}
\]}
\item \textit{Let $a\in \Z$ and $u,v\in \Wordso$ be words such that $u\ta = v\so$. Then in $\Freeoprex$ we have
\[
\begin{array}{r c l}
[u\cdot v]^{\boxtimes a} & = & [u]^{\boxtimes a} \cdot [v]^{\boxtimes a}.
\end{array}
\]}
\item \textit{Let $a\in \Z$ and $w\in \Wordso$ be a word. Then in $\Freeoprex$ we have
\[
\begin{array}{r c l}
[\he_{a,w\so}\op] \cdot [w]^{\boxtimes a} & = & [w] \cdot [\he_{a,w\ta}\op].
\end{array}
\]}
\end{enumerate}
\bpro
Ad 1: This follows by induction over $a$; \cf definition \ref{VER1_DEFI_SE_UND_HE}, definition \ref{VER1_DEFI_POTENZ_IN_WORDS}, definition \ref{VER1_DEFI_FREEOPREX}.

Ad 2: If $a=0$, this follows from definition \ref{VER1_DEFI_SE_UND_HE}. Otherwise using induction over $a$, we conclude that 
\[
\begin{array}{r c l c l}
[u\cdot v]^{\boxtimes a} & \stackrel{\text{\ref{VER1_DEFI_SE_UND_HE}}}{=} & [u\cdot v]^{\boxtimes (a-1)} \boxtimes [u\cdot v]  & \stackrel{\IV}{=} & ([u]^{\boxtimes (a-1)} \cdot [v]^{\boxtimes (a-1)}) \boxtimes ([u] \cdot [v]) \\
 & = & ([u]^{\boxtimes (a-1)} \boxtimes [u]) \cdot ([v]^{\boxtimes (a-1)} \boxtimes [v]) & \stackrel{\text{\ref{VER1_DEFI_SE_UND_HE}}}{=} & [u]^{\boxtimes a} \cdot [v]^{\boxtimes a}.
\end{array}
\]

Ad 3: If $k := |w|\leq 1$, this follows from $\he_{a,w\so}\op \cdot w^{\boxtimes a} \stackrel{4.}{\sim} w \cdot \he_{a,w\ta}$; \cf definition \ref{VER1_DEFI_ELEMENTARE_ÄQUIVALENZEN}. Otherwise let $w=f_0\op \cdot w_1 \cdot f_1\op \cdot {\dots} \cdot w_k \cdot f_k\op $ be the standard decomposition. We define the word $u := f_0\op \cdot w_1 \cdot  f_1\op \cdot {\dots} \cdot w_{k-1} \cdot f_{k-1}\op \in \Wordso$ of length $k-1$ and the word $v := w_k \cdot f_k\op \in \Wordso$ of length $1$. Note that $w = u \cdot v$, $w\so = u\so$, $u\ta = v\so$, $v\ta = w\ta$. Using induction over $k$, we conclude that
\[
\begin{array}{r c l c l c l}
[\he_{a,w\so}\op] \cdot [w]^{\boxtimes a} & = & [\he_{a,u\so}\op] \cdot [u\cdot v]^{\boxtimes a} & \stackrel{\text{2.}}{=} & [\he_{a,u\so}\op] \cdot [u]^{\boxtimes a} \cdot [v]^{\boxtimes a} & \stackrel{\IV}{=} & [u] \cdot [\he_{a,u\ta}] \cdot [v]^{\boxtimes a} \\
 & = & [u] \cdot [\he_{a,v\so}] \cdot [v]^{\boxtimes a} & \stackrel{|v| = 1}{=} & [u]\cdot [v] \cdot [\he_{a,v\ta}\op] & = & [w] \cdot [\he_{a,w\ta}\op].
\end{array}
\]
\epro
\elem

\bdefi \label{VER1_DEFI_STRUKTURMORPHISMUS_FÜR_FREEOX}
We define the morphism of preoperads
\[
\begin{array}{r c c c l}
\mathfrak{f}_{0,X} & : & \mapoop & \to & \Freeoprex \\
 & & f\op & \mapsto & [f\op].
\end{array}
\]
Note that $\mathfrak{f}_{0,X}$ is in fact a morphism of preoperads by remark \ref{VER1_BEM_RECHENREGELN_FÜR_WÖRTER_DER_LÄNGE_NULL}.
\edefi

\blem \label{VER1_LEM_STRUKTURMORPHISMUS_FÜR_FREEOX}
\textit{The morphism of preoperads $\mathfrak{f}_{0,X} : \mapoop \to \Freeoprex$ given by definition \ref{VER1_DEFI_STRUKTURMORPHISMUS_FÜR_FREEOX} satisfies the Braiding Property and the Branching Property; \cf definition \ref{VER1_DEFI_OPERAD}.}
\bpro 

Braiding Property: Let $w,w'\in \Wordso$ be words. Then in $\Freeoprex$ we have
\[
\begin{array}{c l c l}
 & \se_{w\so,w'\so}\op \mathfrak{f}_{0,X} \cdot ([w]\boxtimes [w'])  & = & [\se_{w\so,w'\so}\op \cdot (w\triangleright w'\so) \cdot (w\ta \triangleleft w')] \\
 \stackrel{\ref{VER1_LEM_VORBEREITUNG_SE_STRUKTURMORPHISMUS_VON_FREEOX}}{=} & [(w'\so \triangleleft w) \cdot \se_{w\ta,w'\so}\op\cdot (w\ta \triangleleft w')] & \stackrel{\ref{VER1_LEM_VORBEREITUNG_SE_STRUKTURMORPHISMUS_VON_FREEOX}}{=} & [(w'\so \triangleleft w) \cdot (w' \triangleright w\ta) \cdot \se_{w\ta,w'\ta}\op] \\
 \stackrel{\ref{VER1_LEM_FREEOPREX_INTERCHANGING_RULE}}{=} & [(w' \triangleright w\so) \cdot (w'\ta \triangleleft w) \cdot \se_{w\ta,w'\ta}\op]  & = & ([w']\boxtimes [w]) \cdot \se_{w\ta,w'\ta}\op\mathfrak{f}_{0,X}.
\end{array}
\]

Branching Property: Let $a\in\Z$ and $w\in \Wordso$ be a word. Then in $\Freeoprex$ we have
\[
\begin{array}{r c c c c c l}
\he_{a,w\so}\op \mathfrak{f}_{0,X} \cdot [w]^{\boxtimes a} & = & [\he_{a,w\so}\op] \cdot [w]^{\boxtimes a} & \stackrel{\text{\ref{VER1_LEM_VORBEREITUNG_STRUKTURMORPHISMUS_VON_FREEOX}.3}}{=} & [w] \cdot [\he_{a,w\ta}\op] & = &  [w] \cdot \he_{a,w\ta}\op \mathfrak{f}_{0,X}.
\end{array}
\]
\epro
\elem

\bdefi \label{VER1_DEFI_FREEOX}
We define the \textit{free operad} (over $X$)
\[
\begin{array}{r c l}
\Freeox & := & (\Freeoprex, \mathfrak{f}_{0,X});
\end{array}
\]
\cf definition \ref{VER1_DEFI_FREEOPREX}, definition \ref{VER1_DEFI_STRUKTURMORPHISMUS_FÜR_FREEOX}. By lemma \ref{VER1_LEM_STRUKTURMORPHISMUS_FÜR_FREEOX}, $\Freeox$ is in fact an operad.
\edefi

\subsection{Universal property of $\Freeox$}

\bdefi \label{VER1_DEFI_INKLUSION_VON_X_NACH_FREEOX}
We define the morphism of biindexed sets
\[
\begin{array}{r c c c l}
\biifree_{0,X} & : & X & \to & \Freeoprex \\
 & & x & \mapsto & [\big(0,x,0\big)].
\end{array}
\]
We will often write $[x] \stackrel{\text{short}}{=} x\biifree_{0,X}$.
\edefi

\begin{satz}[Universal property of $\Freeox$] \label{VER1_THEO_UNIV_EIG_VON_FREEOX}
\textit{Let $\oppo$ be an operad and $\biipropp_0 : X \to \proppo$ a morphism of biindexed sets; \cf definition \ref{VER1_DEFI_BIINDEXED_SETS}. Then there exists a unique morphism of operads $\chi : \Freeox \to \oppo$ such that $\biifree_{0,X} \cdot \chi = \biipropp_0$. In particular the following diagram commutes:} 
\[
\begin{tikzcd}[row sep=large, column sep = large]
X \arrow{r}{\biipropp_0} \arrow[swap]{d}{\biifree_{0,X}} & \proppo  \\
\Freeoprex \arrow{ur}{\chi}  & \mapoop \arrow[swap]{u}{\mathfrak{p}_0} \arrow{l}{\mathfrak{f}_{0,X}}
\end{tikzcd}
\]
\bpro Existence: First we construct a functor $\hat{\chi} : \Wordso \to \proppo$ that acts identically on the objects.

Let $w\in \Wordso$ be a word of length $k$ and $w = f_0\op \cdot w_1 \cdot f_1\op \cdot  {\dots} \cdot w_k \cdot f_k\op$ the standard decomposition. Note that there exist unique $l_i,r_i\in \Z$ and $x_i\in X$ for $i\in[1,k]$ such that $w_i = \big(l_i,x_i,r_i\big)$. We define
\[
\begin{array}{r c l}
w\hat{\chi} & = & (f_0\op \cdot \big(l_1,x_1,r_1\big) \cdot f_1\op \cdot  {\dots} \cdot \big(l_k,x_k,r_k\big) \cdot f_k\op) \hat{\chi} \\
 & := & f_0\op\mathfrak{p}_0 \cdot \Big(\id_{l_1}\boxtimes \, x_1\biipropp_0 \boxtimes \id_{r_1}\Big) \cdot f_1\op\mathfrak{p}_0 \cdot {\dots} \cdot \Big(\id_{l_k}\boxtimes \, x_k\biipropp_0 \boxtimes \id_{r_k}\Big)\cdot f_k\op\mathfrak{p}_0.
\end{array}
\]
So in particular we have $f\op\hat{\chi} = f\op\mathfrak{p}_0$ for $f\in \mapo$.
%We define $f_i\op \hat{\chi} := f_i\op\mathfrak{p}_0$ for $i\in [0,k]$ and $w_i\hat{\chi} := \id_{l_i}\boxtimes \, x_i\biipropp_0 \boxtimes \id_{r_i}$ for $i\in[1,k]$. Further we define $w\hat{\chi}$ to be the composition. So in particular we have
%\[
%\begin{array}{r c l}
%w\hat{\chi} & := & f_0\op\mathfrak{p}_0,
%\end{array}
%\]
%if $k=0$, and
%\[
%\begin{array}{r c l}
%w\hat{\chi} & := & f_0\op\mathfrak{p}_0 \cdot \Big(\id_{l_1}\boxtimes \, x_1\biipropp_0 \boxtimes \id_{r_1}\Big) \cdot f_1\op\mathfrak{p}_0 \cdot {\dots} \cdot \Big(\id_{l_k}\boxtimes \, x_k\biipropp_0 \boxtimes \id_{r_k}\Big)\cdot f_k\op\mathfrak{p}_0 ,
%\end{array}
%\]
%if $k\geq 1$. 

Note that $\hat{\chi}$ respects source and target since $\mathfrak{p}_0$ acts identically on the objects. Note that $\hat{\chi}$ respects identities and composition since $\mathfrak{p}_0$ is a functor; \cf the composition in $\Wordso$ given in definition \ref{VER1_DEFI_WORDS}. 

\textit{Claim 1: We have $(q\triangleleft w \triangleright p)\hat{\chi}  =  \id_q \boxtimes \, w\hat{\chi} \boxtimes \id_p$ for $q,p\in \Z$ and a word $w\in \Wordso$.}

This follows from the definition of $q\triangleleft w \triangleright p$, the definition of $\hat{\chi}$ and the fact that $\proppo$ is a preoperad and that $\mathfrak{p}_0$ is a morphism of preoperads. This proves \textit{Claim 1}.

\textit{Claim 2: We have $w^{\boxtimes a}\hat{\chi} = (w\hat{\chi})^{\boxtimes a}$ for $a\in \Z$ and a word $w\in \Wordso$.}

Let $a\in \Z$ and $w\in \Wordso$ be a word. If $a=0$, this follows from definition \ref{VER1_DEFI_POTENZ_IN_WORDS} and definition \ref{VER1_DEFI_SE_UND_HE}. Otherwise using induction over $a$, we conclude that
\[
\begin{array}{c l c l}
 & w^{\boxtimes a}\hat{\chi} & \stackrel{\text{\ref{VER1_DEFI_POTENZ_IN_WORDS}}}{=} & \big((w^{\boxtimes (a-1)}\triangleright w\so) \cdot (w^{\boxtimes (a-1)}\ta \triangleleft w)\big)\hat{\chi} \\
 = & (w^{\boxtimes (a-1)}\triangleright w\so)\hat{\chi} \cdot (w^{\boxtimes (a-1)}\ta \triangleleft w)\hat{\chi}  &  \stackrel{\text{Claim 1}}{=} & (w^{\boxtimes (a-1)}\hat{\chi} \boxtimes \id_{w\so}) \cdot (\id_{w^{\boxtimes (a-1)}\ta} \boxtimes\, w \hat{\chi}) \\
 = & (w^{\boxtimes (a-1)}\hat{\chi} \boxtimes \id_{(w\hat{\chi})\so}) \cdot (\id_{(w^{\boxtimes (a-1)}\hat{\chi})\ta} \boxtimes\, w \hat{\chi})  & \stackrel{\IV}{=} & \big((w\hat{\chi})^{\boxtimes (a-1)} \boxtimes \id_{(w\hat{\chi})\so}\big) \cdot (\id_{(w\hat{\chi})^{\boxtimes (a-1)}\ta} \boxtimes\, w \hat{\chi}) \\
 = & (w\hat{\chi})^{\boxtimes (a-1)} \boxtimes w \hat{\chi}  & \stackrel{\text{\ref{VER1_DEFI_SE_UND_HE}}}{=} & (w\hat{\chi})^{\boxtimes a}.
\end{array}
\]
This proves \textit{Claim 2}.

\textit{Claim 3: We have $w\hat{\chi} = w'\hat{\chi}$ for $w\sim w'$.}

Let $w,w'\in \Wordso$ be words such that $w\sim w'$; \cf definition \ref{VER1_DEFI_ELEMENTARE_ÄQUIVALENZEN}. 

If $w\stackrel{1.}{\sim} w'$, we may choose words $v,v'\in\Wordso$ (of length $\leq 1$) such that $w = (v \triangleright v'\so) \cdot (v\ta \triangleleft v')$ and $w' = (v\so \triangleleft v') \cdot (v \triangleright v'\ta)$. Using the functoriality of $\hat{\chi}$ we have $w\hat{\chi} \stackrel{\text{Claim 1}}{=} v\hat{\chi} \boxtimes v'\hat{\chi} \stackrel{\text{Claim 1}}{=} w'\hat{\chi}$.

If $w\stackrel{2.}{\sim} w'$, we may choose $a,q,p\in \Z$ and a word $v\in\Wordso$ (of length $\leq 1$) such that $w = q \triangleleft \big(\se_{v\so,a}\op \cdot (v \triangleright a)\big) \triangleright p$ and $w' = q \triangleleft \big((a \triangleleft v) \cdot \se_{v\ta,a}\op\big) \triangleright p$. Using the functoriality of $\hat{\chi}$ we have
\[
\begin{array}{r c c c c c l}
w\hat{\chi} \hspace{-1ex} & \hspace{-1ex} \stackrel{\text{Claim 1}}{=} \hspace{-1ex} &  \hspace{-1ex} \id_q \boxtimes \, \big(\se_{v\so,a}\op\mathfrak{p}_0 \cdot (v\hat{\chi} \boxtimes \id_a)\big) \boxtimes \id_p \hspace{-1ex} & \hspace{-1ex} \stackrel{\text{Braiding}}{=} \hspace{-1ex} & \hspace{-1ex} \id_q \boxtimes \, \big((\id_a \boxtimes \, v\hat{\chi}) \cdot \se_{v\ta,a}\op\mathfrak{p}_0\big) \boxtimes \id_p \hspace{-1ex} & \hspace{-1ex} \stackrel{\text{Claim 1}}{=} \hspace{-1ex} & \hspace{-1ex} w'\hat{\chi}.
\end{array}
\]

If $w\stackrel{3.}{\sim} w'$, we may choose $a,q,p\in \Z$ and a word $v\in\Wordso$ (of length $\leq 1$) such that $w = q \triangleleft \big(\se_{a,v\so}\op \cdot (a \triangleleft v)\big) \triangleright p$ and $w' = q \triangleleft \big((v \triangleright a) \cdot \se_{a,v\ta}\op\big) \triangleright p$. Using the functoriality of $\hat{\chi}$ we have 
\[
\begin{array}{r c c c c c l}
w\hat{\chi} \hspace{-1ex} & \hspace{-1ex} \stackrel{\text{Claim 1}}{=} \hspace{-1ex} & \hspace{-1ex} \id_q \boxtimes \, \big(\se_{a,v\so}\op\mathfrak{p}_0 \cdot (\id_a \boxtimes \, v\hat{\chi})\big) \boxtimes \id_p \hspace{-1ex} & \hspace{-1ex} \stackrel{\text{Braiding}}{=} \hspace{-1ex} & \hspace{-1ex} \id_q \boxtimes \, \big((v\hat{\chi} \boxtimes \id_a ) \cdot \se_{a,v\ta}\op\mathfrak{p}_0\big) \boxtimes \id_p \hspace{-1ex} & \hspace{-1ex} \stackrel{\text{Claim 1}}{=} \hspace{-1ex} & \hspace{-1ex} w'\hat{\chi}.
\end{array}
\]

If $w\stackrel{4.}{\sim} w'$, we may choose $a,q,p\in \Z$ and a word $v\in\Wordso$ (of length $\leq 1$) such that $w = q \triangleleft (\he_{a,v\so}\op \cdot v^{\boxtimes a}) \triangleright p$ and $w' = q \triangleleft (v \cdot \he_{a,v\ta}\op) \triangleright p$. Using the functoriality of $\hat{\chi}$ we have 
\[
\begin{array}{c l c l}
 & w\hat{\chi} & \stackrel{\text{Claim 1}}{=} &  \id_q \boxtimes \, (\he_{a,v\so}\op\mathfrak{p}_0 \cdot v^{\boxtimes a}\hat{\chi}) \boxtimes \id_p \\
 \stackrel{\text{Claim 2}}{=} & \id_q \boxtimes \, \big(\he_{a,v\so}\op\mathfrak{p}_0 \cdot (v\hat{\chi})^{\boxtimes a}\big) \boxtimes \id_p & \stackrel{\text{Branching}}{=} & \id_q \boxtimes \, (v\hat{\chi} \cdot \he_{a,v\ta}\op\mathfrak{p}_0) \boxtimes \id_p \\
  \stackrel{\text{Claim 1}}{=} & w'\hat{\chi}. 
\end{array}
\]

This proves \textit{Claim 3}.

Claim 3 yields the functor
\[
\begin{array}{r c c c l}
\chi & : & \Freeoprex & \to & \proppo \\
 & & [w] & \mapsto & w\hat{\chi}
\end{array}
\]
by the universal property of the factor category. Note that $\chi$ acts identically on the objects.

\textit{Claim 4: The functor $\chi$ is a morphism of preoperads.}

Let $w,w'\in \Wordso$ be words. Then we have 
\[
\begin{array}{c l c l}
 & ([w]\boxtimes [w'])\chi & = & [(w\triangleright w'\so)\cdot (w\ta \triangleleft w')]\chi \\
 = & \big((w\triangleright w'\so)\cdot (w\ta \triangleleft w')\big)\hat{\chi} & = & (w\triangleright w'\so)\hat{\chi}\cdot (w\ta \triangleleft w')\hat{\chi} \\
 \stackrel{\text{Claim 1}}{=} & (w\hat{\chi}\boxtimes \id_{w'\so})\cdot (\id_{w\ta} \boxtimes\, w'\hat{\chi}) & = &  w\hat{\chi}\boxtimes  w'\hat{\chi} \\
 = & [w]\chi\boxtimes  [w']\chi.
\end{array}
\]
This proves \textit{Claim 4}. 

\textit{Claim 5: We have $\mathfrak{f}_{0,X} \cdot \chi = \mathfrak{p}_0$.} 

Let $f\in \mapo$. Then we have $f\op(\mathfrak{f}_{0,X} \cdot \chi) = [f\op]\chi  =  f\op\hat{\chi}  =  f\op\mathfrak{p}_0$. This proves \textit{Claim 5}.

By Claim 4 and Claim 5, $\chi : \Freeox  \to  \oppo$ is a morphism of operads.

\textit{Claim 6: We have $\biifree_{0,X} \cdot \chi = \biipropp_0$.}

Let $x\in X$. Then we have $x(\biifree_{0,X} \cdot \chi)  =  [\big(0,x,0\big)]\chi  =  \big(0,x,0\big) \hat{\chi}  =  \id_0\boxtimes \, x\biipropp_0 \boxtimes \id_0  =  x\biipropp_0$; \cf definition \ref{VER1_DEFI_INKLUSION_VON_X_NACH_FREEOX}. This proves \textit{Claim 6}. 

So Claim 6 shows that the morphism of operads $\chi : \Freeox \to \oppo$ constructed above satisfies the properties stated in the theorem.

Uniqueness: Let $l,r\in \Z$ and $x\in X$. Note that in $\Freeoprex$ we have 
\[
\begin{array}{r c c c l}
[\id_l\op]\boxtimes x\biifree_{0,X} \boxtimes [\id_r\op] & = & [\id_l\op]\boxtimes [\big(0,x,0\big)]\boxtimes [\id_r\op] & \stackrel{\text{\ref{VER1_BEM_OXO_MIT_IDENTITÄTEN}}}{=} &  [\big(l,x,r\big)].
\end{array}
\]
Let $\varphi_0 : \Freeox \to \oppo$ be a morphism of operads such that $\biifree_{0,X} \cdot \varphi_0 = \biipropp_0$. Let $w\in\Words$ be a word and $w=f_0\op\cdot w_1 \cdot f_1\op \cdot {\dots} \cdot w_k \cdot f_k\op$ be the standard decomposition. 

Since $\varphi_0$ is a functor we have $[w]\varphi_0  =  [f_0\op]\varphi_0 \cdot [w_1]\varphi_0 \cdot [f_1\op]\varphi_0 \cdot {\dots} \cdot [w_k]\varphi_0 \cdot [f_k\op]\varphi_0$.

Since $\varphi_0$ is a morphism of operads we have $[f_i\op]\varphi_0  =  f_i\op(\mathfrak{f}_{0,X} \cdot \varphi_0)  =   f_i\op\mathfrak{p}_0$ for $i\in[0,k]$.

Note that there exist unique $l_i,r_i\in \Z$ and $x_i\in X$ such that $w_i = \big(l_i,x_i,r_i\big)$. Since $\biifree_{0,X} \cdot \varphi_0 = \biipropp_0$, we conclude that
\[
\begin{array}{r c c c c c l}
[w_i]\varphi_0 & = & [\big(l_i,x_i,r_i\big)]\varphi_0 & \stackrel{\text{\ref{VER1_BEM_OXO_MIT_IDENTITÄTEN}}}{=} & \big([\id_l\op]\boxtimes  x_i\biifree_{0,X} \boxtimes [\id_r\op]\big)\varphi_0 & = & \id_{l_i} \boxtimes \, x_i\biipropp_0 \boxtimes \id_{r_i} 
\end{array}
\]
for $i\in[1,k]$. This proves that $\varphi_0$ equals the morphism of operads $\chi$ constructed above.
\epro
\end{satz}

\section{Generators and relations}

\bdefi \label{VER1_DEFI_FACTOR_PREOPERAD_AND_MORPHISM}
Let $\proppo$ be a preoperad. We define a \textit{multiplicative congruence relation} (on $\proppo$) as a congruence relation $(\equiv) \subseteq \proppo \btimes \proppo$ such that $\id_l\boxtimes \, a \boxtimes \id_r \equiv \id_l\boxtimes \, b \boxtimes \id_r$ for $l,r\in \Z$ and $a\equiv b$; \cf definition \ref{VER1_DEFI_BIINDEXED_SETS}. Note that $x\boxtimes x' \equiv y\boxtimes y'$ for $x,x',y,y'\in \proppo$ such that $x\equiv y$ and $x'\equiv y'$.

Recall that the factor category $\proppo / (\equiv)$ is again a category. We endow $\proppo/(\equiv)$ with the monoidal product
\[
\begin{array}{r c l}
[x] \boxtimes [x'] & := & [x\boxtimes x']
\end{array}
\]
for $x\in \proppo(m,n)$ and $x'\in\proppo(m',n')$, where $m,n,m',n'\in\Z$. Then $\proppo / (\equiv)$ is a preoperad; \cf \bfcite{Eggert}{Def.\ 2.38}.

We define the \textit{factor morphism}
\[
\begin{array}{r c c c l}
\funcfacfunc_{\proppo,(\equiv)} & : & \proppo & \to & \proppo/(\equiv) \\
 & & x & \mapsto & [x].
\end{array}
\]
Note that $\funcfacfunc_{\proppo,(\equiv)}$ is in fact a morphism of preoperads; \cf \bfcite{Eggert}{Def.\ 2.39}.
\edefi

\bdefi \label{VER1_DEFI_EQUIV_Y}
Let $\proppo$ be a preoperad. %We define a \textit{congruence-generating subset} as a subset $Y\subseteq \proppo \times \proppo$ such that $a\so = b\so$ and $a\ta = b\ta$ for $(a,b) \in Y$.
Let $Y\subseteq \proppo \btimes \proppo$ be a subset; \cf definition \ref{VER1_DEFI_BIINDEXED_SETS}. We define $(\equiv_Y) \subseteq \proppo \btimes \proppo$ as the congruence relation generated by 
\[
\begin{array}{r c l}
\hat{Y} & := & \{(\id_l\boxtimes \, a \boxtimes \id_r,\id_l\boxtimes \, b \boxtimes \id_r) \in \proppo \btimes \proppo \,|\, l,r\in \Z, \big((a,b)\in Y \vee (b,a)\in Y\big) \}.
\end{array}
\]
The congruence relation $(\equiv_Y)$ is in fact a multiplicative congruence relation containing $Y$. Moreover for every multiplicative congruence relation $(\diamond) \subset \proppo \btimes \proppo$ containing $Y$ we have $(\equiv_Y) \subset (\diamond)$.
\edefi

\bdefi \label{VER1_DEFI_FACTOR_OPERAD}
Let $\oppo$ be an operad and $(\equiv)\subseteq \proppo \btimes \proppo$ a multiplicative congruence relation. We define the morphism of preoperads 
\[
\begin{array}{r c c c c c l}
\overline{\mathfrak{p}_0} & := & \mathfrak{p}_0 \cdot \funcfacfunc_{\proppo,(\equiv)} & : & \mapoop & \to & \proppo/(\equiv) \\
 & & & & f\op & \mapsto & [f\op\mathfrak{p}_0]. 
\end{array}
\] 
Note that $\overline{\mathfrak{p}_0}$ satisfies the Braiding Property and the Branching Property; \cf definition \ref{VER1_DEFI_OPERAD}. We define the \textit{factor operad}
\[
\begin{array}{r c l}
\oppo / (\equiv) & := & \big(\proppo/(\equiv), \overline{\mathfrak{p}_0}\big).
\end{array}
\]
Recall that we have the factor morphism $\funcfacfunc_{\proppo,(\equiv)}  :  \proppo  \to  \proppo/(\equiv)$. Since $\mathfrak{p}_0 \cdot \funcfacfunc_{\proppo,(\equiv)} = \overline{\mathfrak{p}_0}$, the morphism $\funcfacfunc_{\proppo,(\equiv)} : \oppo \to \oppo/(\equiv)$ is in fact a morphism of operads.
\edefi

\begin{lem}[Universal property of the factor operad] \label{VER1_LEM_UNIV_EIG_FACTOR_OPERAD}
\textit{Let $\oppo$ and $\qoppo$ be operads and $Y\subseteq \proppo \btimes \proppo$ a subset; \cf definition \ref{VER1_DEFI_BIINDEXED_SETS}. Let $\varphi_0 : \oppo \to \qoppo$ be a morphism of operads such that $a\varphi_0  = b\varphi_0$ for $(a,b)\in Y$. Then there exists a unique morphism of operads $\psi_0 : \proppo/(\equiv_Y) \to \qproppo$ such that $\funcfacfunc_{\proppo,(\equiv_Y)} \cdot \, \psi_0 = \varphi_0$.}
\bpro
Since $\varphi_0$ is a morphism of preoperads we conclude that $u\varphi_0 = v\varphi_0$ for $(u,v)\in \hat{Y}$. By the universal property of the factor category, there exists a unique functor $\psi_0 : \proppo/(\equiv_Y) \to \qproppo$ such that $\funcfacfunc_{\proppo,(\equiv_Y)} \cdot \, \psi_0 = \varphi_0$; \cf definition \ref{VER1_DEFI_EQUIV_Y}. Note that $\psi_0$ acts identically on the objects.

Let $x,x'\in \proppo$. Then we have $([x]\boxtimes [x'])\psi_0 = [x\boxtimes x']\psi_0 = (x\boxtimes x')\varphi_0 = x\varphi_0 \boxtimes x'\varphi_0 = [x]\psi_0 \boxtimes [x']\psi_0$. So $\psi_0$ is a morphism of preoperads.

Since $\overline{\mathfrak{p}_0} \cdot \psi_0 \stackrel{\ref{VER1_DEFI_FACTOR_OPERAD}}{=} \mathfrak{p}_0 \cdot \funcfacfunc_{\proppo,(\equiv_Y)} \cdot \, \psi_0 = \mathfrak{p}_0 \cdot \varphi_0 = \mathfrak{q}_0$, we conclude that $\psi_0$ is a morphism of operads.
\epro
\end{lem}

\bdefi \label{VER1_DEFI_OPERADE_X_MODULO_Y}
Let $X$ be a biindexed set. Let $Y\subseteq \Freeoprex \btimes \Freeoprex$ be a subset; \cf definition \ref{VER1_DEFI_BIINDEXED_SETS}. We write
\[
\begin{array}{r c l}
\generated\xy  & := & \Freeox / (\equiv_Y)
\end{array}
\]
for the factor operad; \cf definition \ref{VER1_DEFI_FACTOR_OPERAD}.
\edefi

\section{The Eggert-operad $\Groupfp$}

Throughout this section let $X_{\operatorname{G}} = \big(X_{\operatorname{G}}(m,n)\big)_{m,n\in \Z}$ be the biindexed set given by
\[
\begin{array}{r c c l}
X_{\operatorname{G}}(m,n) & := & \begin{cases}
\{\mu\} & \text{ if } m=2, n=1\\
\{\eta\} & \text{ if } m=0, n=1\\
\{\omega\} & \text{ if } m=1, n=1\\
\emptyset & \text{ else}.
\end{cases}
\end{array}
\] 
Recall that we will use the abbreviation $[x] \stackrel{\text{short}}{=} [\big(0,x,0\big)]  = x\biifree_{0,X_{\operatorname{G}}} \in \Freeopreg$ for $x\in X_{\operatorname{G}}$; \cf definition \ref{VER1_DEFI_INKLUSION_VON_X_NACH_FREEOX}.

We will denote a group as $M = (M,\mu_M, \eta_M,\omega_M)$, where $\mu_M : M^{\times 2} \to M^{\times 1}$ is the group multiplication, $\eta_M : M^{\times 0} \to M^{\times 1}$ is the unit and $\omega_M : M^{\times 1} \to M^{\times 1}$ is the inversion map. So $(m,m')\mu_M = mm'$ is the product, $()\eta_M = 1_M$ is the neutral element and $m\omega_M = m^{-1}$ is the inverse element of $m\in M$.

\subsection{Construction of $\Groupfp$}

We will make use of two non-bijective elements of $\mapo$; \cf remark \ref{VER1_BEM_VERGLEICH_PROP_UND_EGGERT_OPERADE}.

\bdefi \label{VER1_DEFI_FTWO_FZERO_DIAG_CONST}
We will denote the unique map in $\mapo(2,1)$ as $\ftwo$; \cf example \ref{VER1_BEI_MAPO}. Further we will denote the unique map in $\mapo(0,1)$ as $\fzero$.

Let $M$ be a set. Recall that we have the operad $\Eno(M) = (\eno(M),\mathfrak{e}_0)$; \cf example \ref{VER1_BEI_OPERADE_ENO(M)}. We define the maps
\[
\begin{array}{r c c c c c l}
\diag_M & := & {\ftwo}\op\mathfrak{e}_0  & : & M^{\times 1} & \to & M^{\times 2} \\
 & & & & m & \mapsto & (m,m)
\end{array}
\]
and
\[
\begin{array}{r c c c c c l}
\konst_M & := & {\fzero}\op\mathfrak{e}_0  & : & M^{\times 1} & \to & M^{\times 0} \\
 & & & & m & \mapsto & ();
\end{array}
\]
\edefi

\bdefi \label{VER1_DEFI_XG_UND_YG}
We define
\[
\begin{array}{r c l c l}
\vsp (y_{\operatorname{G},1}^1,y_{\operatorname{G},2}^1) & := & \big(([\mu]\boxtimes [\id_1\op])\cdot [\mu] , ([\id_1\op]\boxtimes  [\mu])\cdot [\mu]\big) & \in & \big(\Freeopreg(3,1)\big)^{\times 2} \\
\vsp (y_{\operatorname{G},1}^2,y_{\operatorname{G},2}^2) & := & \big(([\eta]\boxtimes [\id_1\op]) \cdot [\mu] , [\id_1\op]\big) & \in & \big(\Freeopreg(1,1)\big)^{\times 2} \\
\vsp (y_{\operatorname{G},1}^3,y_{\operatorname{G},2}^3) & := & \big(([\id_1\op]\boxtimes  [\eta]) \cdot [\mu] , [\id_1\op]\big) & \in & \big(\Freeopreg(1,1)\big)^{\times 2} \\
\vsp (y_{\operatorname{G},1}^4,y_{\operatorname{G},2}^4) & := & \big([{\ftwo}\op] \cdot ([\omega] \boxtimes [\id_1\op]) \cdot [\mu] , [{\fzero}\op] \cdot [\eta] \big) & \in & \big(\Freeopreg(1,1)\big)^{\times 2} \\
(y_{\operatorname{G},1}^5,y_{\operatorname{G},2}^5) & := & \big( [{\ftwo}\op] \cdot ([\id_1\op] \boxtimes  [\omega]) \cdot [\mu] , [{\fzero}\op] \cdot [\eta] \big) & \in & \big(\Freeopreg(1,1)\big)^{\times 2}.
\end{array}
\]
Note that $\id_{m,\Freeoprex}  =  [\id_m\op]  \in  \Freeopreg(m,m)$ for $m\in \Z$. Further we define
\[
\begin{array}{r c c c l}
Y_{\operatorname{G}} & := & \{(y_{\operatorname{G},1}^i,y_{\operatorname{G},2}^i)\,|\, i\in[1,5]\} & \subseteq & \Freeopreg \btimes \Freeopreg.
\end{array}
\]
We define the operad $\Groupfp$ as 
\[
\begin{array}{r c c c l}
\Groupfp & := & \big(\Groupfppre, \mathfrak{g}_0\fp\big) & := & \generated\xgyg;
\end{array}
\]
\cf definition \ref{VER1_DEFI_OPERADE_X_MODULO_Y}. We add the upper index ``fp'' to indicate that $\Groupfp$ is finitely presented, $\ie$ that both $X_{\operatorname{G}}$ and $Y_{\operatorname{G}}$ are finite. In accordance with the notation of Eggert, we add an upper index ``pre'' to indicate the underlying preoperad of an operad. Recall that we have the factor morphism
\[
\begin{array}{r c c c l}
\funcfacfunc_{\Freeopreg,(\equiv_{Y_{\operatorname{G}}})} & : & \Freeog & \to & \Groupfp\\
 & & [w] & \mapsto & [[w]];
\end{array}
\]
\cf definition \ref{VER1_DEFI_FACTOR_PREOPERAD_AND_MORPHISM}.
\edefi

\bsatz 
\textit{Let $\big(M,\gamma\big)$ be a $\Groupfp$-algebra; \cf definition \ref{VER1_DEFI_ALGEBRA}. So we have the morphism of operads}
\[
\begin{array}{r c c c c c l}
\gamma & : & \Groupfp & \to & \Eno(M) & & \\
 & & [[\mu]] & \mapsto & [[\mu]]\gamma & \in & \eno(M)(2,1)\\
 & & [[\eta]] & \mapsto & [[\eta]]\gamma & \in & \eno(M)(0,1)\\
 & & [[\omega]] & \mapsto & [[\omega]]\gamma & \in & \eno(M)(1,1).
\end{array}
\]
\textit{Let}
\[
\begin{array}{r c c c l c l}
\mu_M & := &  [[\mu]]\gamma & : & M^{\times 2} & \to & M^{\times 1} \\
\eta_M & := &  [[\eta]]\gamma & : & M^{\times 0} & \to & M^{\times 1} \\
\omega_M & := &  [[\omega]]\gamma & : & M^{\times 1} & \to & M^{\times 1}. 
\end{array}
\]
\textit{Then $\big(M,\mu_M, \eta_M,\omega_M\big)$ is a group.}
\bpro
Since $\gamma$ is a morphism of operads, we have
\[
\begin{array}{r c c c c c l}
\diag_M & = & {\ftwo}\op\mathfrak{e}_0 & = & {\ftwo}\op(\mathfrak{g}_0\fp \cdot \gamma) & = & [[{\ftwo}\op]]\gamma
\end{array}
\]
and
\[
\begin{array}{r c c c c c l}
\konst_M & = & {\fzero}\op\mathfrak{e}_0 & = & {\fzero}\op(\mathfrak{g}_0\fp \cdot \gamma) & = & [[{\fzero}\op]]\gamma.
\end{array}
\]
Recall that $[[\id_1\op]]\gamma = \id_{1,\eno(M)} = \id_M : M \to M$. We obtain the following equations in $\eno(M)$. 

$(\mu_M\boxtimes \id_M)\cdot \mu_M = \big(([[\mu]]\boxtimes [[\id_1\op]])\cdot[[\mu]]\big)\gamma \stackrel{\text{\ref{VER1_DEFI_XG_UND_YG}}}{=} \big(([[\id_1\op]]\boxtimes [[\mu]])\cdot[[\mu]]\big)\gamma = (\id_M\boxtimes\, \mu_M)\cdot \mu_M$

So $\mu_M$ is associative.

$(\eta_M\boxtimes \id_M)\cdot \mu_M = \big(([[\eta]] \boxtimes [[\id_1\op]])\cdot[[\mu]]\big)\gamma \stackrel{\text{\ref{VER1_DEFI_XG_UND_YG}}}{=} [[\id_1\op]]\gamma =  \id_M$

So $\eta_M$ is left-neutral.

$(\id_M\boxtimes \, \eta_M)\cdot \mu_M = \big(([[\id_1\op]] \boxtimes  [[\eta]])\cdot[[\mu]]\big)\gamma \stackrel{\text{\ref{VER1_DEFI_XG_UND_YG}}}{=} [[\id_1\op]]\gamma = \id_M$

So $\eta_M$ is right-neutral.

$\diag_M \cdot (\omega_M \boxtimes \id_M) \cdot \mu_M = \big([[{\ftwo}\op]] \cdot ([[\omega]] \boxtimes [[\id_1\op]]) \cdot [[\mu]]\big)\gamma \stackrel{\text{\ref{VER1_DEFI_XG_UND_YG}}}{=} \big([[{\fzero}\op]] \cdot [[\eta]]\big)\gamma = \konst_M \cdot \eta_M$

So $\omega_M$ is left-inverse.

$\diag_M \cdot (\id_M \boxtimes \, \omega_M) \cdot \mu_M = \big([[{\ftwo}\op]] \cdot ([[\id_1\op]] \boxtimes  [[\omega]]) \cdot [[\mu]]\big)\gamma \stackrel{\text{\ref{VER1_DEFI_XG_UND_YG}}}{=} \big([[{\fzero}\op]] \cdot [[\eta]]\big)\gamma = \konst_M \cdot \eta_M$

So $\omega_M$ is right-inverse. We conclude that $\big(M,\mu_M, \eta_M,\omega_M\big)$ is a group.
\epro
\esatz

\bsatz 
\textit{Let $\big(M,\mu_M, \eta_M,\omega_M\big)$ be a group. Then we have the morphism of operads
\[
\begin{array}{r c c c l}
\gamma & : & \Groupfp & \to & \Eno(M) \\
 & & [[\mu]] & \mapsto & \mu_M \\
 & & [[\eta]] & \mapsto & \eta_M \\
 & & [[\omega]] & \mapsto & \omega_M.
\end{array}
\]
So $\big(M,\gamma\big)$ is a $\Groupfp$-algebra.}
\bpro
We define the morphism of biindexed sets
\[
\begin{array}{r c c c l}
\gamma_{\mathrm{bii}} & : &  X_{\operatorname{G}} & \to & \eno(M)\\
 & & \mu & \mapsto & \mu_M \\
 & & \eta & \mapsto & \eta_M \\
 & & \omega & \mapsto & \omega_M.
\end{array}
\]
By the universal property of $\biifree_{0,X_{\operatorname{G}}} : X_{\operatorname{G}} \to \Freeopreg$, there exists a unique morphism of operads $\gamma_\mathrm{fr}  :  \Freeog  \to  \Eno(M)$ such that $\biifree_{0,X_{\operatorname{G}}} \cdot \gamma_\mathrm{fr}  =  \gamma_{\mathrm{bii}}$, $\ie$
\[
\begin{array}{r c c c l}
\gamma_\mathrm{fr} & : &  \Freeog & \to & \Eno(M) \\
 & & [\mu] & \mapsto & \mu_M \\
 & & [\eta] & \mapsto & \eta_M \\
 & & [\omega] & \mapsto & \omega_M;
\end{array}
\]
\cf theorem \ref{VER1_THEO_UNIV_EIG_VON_FREEOX}.
%$\ie$ such that 
%\[
%\begin{array}{r c l}
%\biifree_{0,X_{\operatorname{G}}} \cdot \gamma_\mathrm{fr} & = & \gamma_{\mathrm{bii}}.
%\end{array}
%\]
Since $\gamma_\mathrm{fr}$ is a morphism of operads, we have
\[
\begin{array}{r c c c c c l}
\diag_M  & = & {\ftwo}\op\mathfrak{e}_0 & = & {\ftwo}\op(\mathfrak{f}_{0,X_{\operatorname{G}}} \cdot \gamma_\mathrm{fr}) & \stackrel{\ref{VER1_DEFI_STRUKTURMORPHISMUS_FÜR_FREEOX}}{=} & [{\ftwo}\op]\gamma_\mathrm{fr}
\end{array}
\]
and
\[
\begin{array}{r c c c c c l}
\konst_M  & = & {\fzero}\op\mathfrak{e}_0 & = & {\fzero}\op(\mathfrak{f}_{0,X_{\operatorname{G}}} \cdot \gamma_\mathrm{fr}) & \stackrel{\ref{VER1_DEFI_STRUKTURMORPHISMUS_FÜR_FREEOX}}{=} & [{\fzero}\op]\gamma_\mathrm{fr}.
\end{array}
\]
We show that 
\[
\begin{array}{r c c c l}
y_{\operatorname{G},1}^i\gamma_\mathrm{fr} & = & y_{\operatorname{G},2}^i\gamma_\mathrm{fr} & \in & \eno(M)
\end{array}
\]
for $i\in [1,5]$; \cf definition \ref{VER1_DEFI_XG_UND_YG}. We obtain the following equations in $\eno(M)$.

Since $\mu_M$ is associative we conclude that

$\big(([\mu]\boxtimes [\id_1\op])\cdot [\mu]\big)\gamma_\mathrm{fr} = (\mu_M\boxtimes \id_M)\cdot \mu_M  = (\id_M\boxtimes \, \mu_M)\cdot \mu_M =  \big(([\id_1\op]\boxtimes [\mu])\cdot [\mu]\big)\gamma_\mathrm{fr}$.

Since $\eta_M$ is left-neutral we conclude that

$\big(([\eta] \boxtimes [\id_1\op])\cdot [\mu]\big)\gamma_\mathrm{fr} =(\eta_M\boxtimes \id_M)\cdot \mu_M  =  \id_M =  [\id_1\op]\gamma_\mathrm{fr}$.

Since $\eta_M$ is right-neutral we conclude that

$\big(([\id_1\op] \boxtimes [\eta])\cdot [\mu]\big)\gamma_\mathrm{fr} = (\id_M\boxtimes \,\eta_M)\cdot \mu_M = \id_M = [\id_1\op]\gamma_\mathrm{fr}$.

Since $\omega_M$ is left-inverse we conclude that

$\big([{\ftwo}\op] \cdot ([\omega] \boxtimes [\id_1\op]) \cdot [\mu]\big)\gamma_\mathrm{fr} = \diag_M \cdot (\omega_M \boxtimes \id_M) \cdot \mu_M = \konst_M \cdot \eta_M = \big([{\fzero}\op] \cdot [\eta]\big)\gamma_\mathrm{fr}$.

Since $\omega_M$ is right-inverse we conclude that

$\big([{\ftwo}\op] \cdot ([\id_1\op] \boxtimes  [\omega]) \cdot [\mu]\big)\gamma_\mathrm{fr} = \diag_M \cdot (\id_M \boxtimes \, \omega_M) \cdot \mu_M = \konst_M \cdot \eta_M = \big([{\fzero}\op] \cdot [\eta]\big)\gamma_\mathrm{fr}$.

So we have 
\[
\begin{array}{r c c c l}
y_{\operatorname{G},1}^i\gamma_\mathrm{fr} & = & y_{\operatorname{G},2}^i\gamma_\mathrm{fr} & \in & \eno(M)
\end{array}
\]
for $i\in [1,5]$. By the universal property of $\funcfacfunc_{\Freeopreg,(\equiv_{Y_{\operatorname{G}}})}  :  \Freeog  \to  \Groupfp$, there exists a unique morphism of operads $\gamma  :  \Groupfp  \to  \Eno(M)$ such that $\funcfacfunc_{\Freeopreg,(\equiv_{Y_{\operatorname{G}}})} \cdot \, \gamma  =  \gamma_\mathrm{fr}$, $\ie$
\[
\begin{array}{r c c c l}
\gamma & : & \Groupfp & \to & \Eno(M)\\
 & & [[\mu]] & \mapsto & \mu_M \\
 & & [[\eta]] & \mapsto & \eta_M \\
 & & [[\omega]] & \mapsto & \omega_M;
\end{array}
\]
\cf lemma \ref{VER1_LEM_UNIV_EIG_FACTOR_OPERAD}.
%$\ie$ such that 
%\[
%\begin{array}{r c l}
%\funcfacfunc_{\Freeopreg,(\equiv_{Y_{\operatorname{G}}})} \cdot \, \gamma & = & \gamma_\mathrm{fr}.
%\end{array}
%\]
So $\big(M,\gamma\big)$ is a $\Groupfp$-algebra.
\epro
\esatz

\subsection{Properties of $\Groupfp$}

\blem \label{VER1_LEM_LINKSAXIOME_BEI_GRUPPEN_REICHEN}
\textit{Recall that we have the unique maps $\fzero, \ftwo\in \mapo$; \cf definition \ref{VER1_DEFI_FTWO_FZERO_DIAG_CONST}. The following equations hold in $\Freeopreg$.}
\begin{enumerate}
\item \hfil $[{\ftwo}\op ]\cdot \big([{\ftwo}\op]\boxtimes [\id_1\op]\big) \;\;\,  = \;\;\, [{\ftwo}\op] \cdot \big([\id_1\op] \boxtimes [{\ftwo}\op]\big)$ \hfill
\item \hfil $[{\ftwo}\op] \cdot \big([{\ftwo}\op]\boxtimes  [{\ftwo}\op]\big) \;\;\, = \;\;\, [{\ftwo}\op] \cdot \big([{\ftwo}\op] \boxtimes  [\id_1\op]\big) \cdot \big([\id_1\op] \boxtimes  [{\ftwo}\op]  \boxtimes  [\id_1\op]\big)$ \hfill
\item \hfil $[{\ftwo}\op] \cdot \big([{\fzero}\op]\boxtimes [\id_1\op]\big) \;\;\, = \;\;\, [\id_1\op] \;\;\, = \;\;\, [{\ftwo}\op] \cdot \big([\id_1\op] \boxtimes [{\fzero}\op]\big)$ \hfill
\item \hfil $[{\ftwo}\op] \cdot ([\omega]\boxtimes [\omega]) \;\;\, = \;\;\, [\omega] \cdot [{\ftwo}\op]$ \hfill
\item \hfil $[\omega] \cdot [{\fzero}\op]  \;\;\, = \;\;\, [{\fzero}\op]$ \hfill
\end{enumerate}
Mapping these equations to $\Groupfp$ via $\funcfacfunc_{\Freeopreg,(\equiv_Y)}$, we obtain according equations in $\Groupfp$. 
\bpro Ad 1: Note that $(\ftwo\boxtimes \id_1) \cdot \ftwo  =  (\id_1\boxtimes\, \ftwo) \cdot \ftwo  \in  \mapo(3,1)$ is the unique map. We conclude that in $\mapoop$ we have 
\[
\begin{array}{c l c l c l}
 & {\ftwo}\op \cdot \big({\ftwo}\op \boxtimes \id_1\op\big) & = & {\ftwo}\op \cdot (\ftwo\boxtimes \id_1)\op  & = & \big((\ftwo\boxtimes \id_1) \cdot \ftwo\big)\op \\
 = & \big((\id_1\boxtimes\, \ftwo) \cdot \ftwo\big)\op  & = & {\ftwo}\op \cdot (\id_1 \boxtimes\, \ftwo)\op  & = & {\ftwo}\op \cdot \big(\id_1\op \boxtimes\, {\ftwo}\op\big).
\end{array}
\]

Ad 2: Note that $(\ftwo \boxtimes \ftwo) \cdot \ftwo  =  (\id_1 \boxtimes\, \ftwo \boxtimes \id_1) \cdot (\ftwo \boxtimes \id_1) \cdot \ftwo  \in  \mapo(4,1)$ is the unique map. We conclude that in $\mapoop$ we have
\[
\begin{array}{c l c l}
 & {\ftwo}\op \cdot \big({\ftwo}\op\boxtimes {\ftwo}\op\big) & = & {\ftwo}\op \cdot (\ftwo\boxtimes \ftwo)\op \\
 = & \big((\ftwo \boxtimes \ftwo) \cdot \ftwo\big)\op  & = & \big((\id_1 \boxtimes\, \ftwo \boxtimes \id_1) \cdot (\ftwo \boxtimes \id_1) \cdot \ftwo\big)\op \\
 = & {\ftwo}\op \cdot (\ftwo \boxtimes \id_1)\op \cdot (\id_1 \boxtimes\, \ftwo \boxtimes \id_1)\op  & = & {\ftwo}\op \cdot \big({\ftwo}\op \boxtimes \id_1\op\big) \cdot \big(\id_1\op \boxtimes\, {\ftwo}\op \boxtimes \id_1\op\big).
\end{array}
\] 

Ad 3: Note that $\id_1  =  (\fzero\boxtimes \id_1) \cdot \ftwo  =  (\id_1\boxtimes\, \fzero) \cdot \ftwo  \in  \mapo(1,1)$ is the unique map. We conclude that in $\mapoop$ we have
\[
\begin{array}{c l c l c l c l}
 & {\ftwo}\op \cdot \big({\fzero}\op \boxtimes \id_1\op\big) & = & {\ftwo}\op \cdot (\fzero\boxtimes \id_1)\op  & = & \big((\fzero\boxtimes \id_1) \cdot \ftwo\big)\op  & = & \id_1\op  \\
 = & \big((\id_1\boxtimes\, \fzero) \cdot \ftwo\big)\op  & = & {\ftwo}\op \cdot (\id_1 \boxtimes\, \fzero)\op  & = & {\ftwo}\op \cdot \big(\id_1\op \boxtimes\, {\fzero}\op\big).
\end{array}
\]

Ad 4: Note that $\ftwo  =  \he_{2,1}  \in  \mapo(2,1)$ is the unique map. In $\Freeoprex$, we obtain
\[
\begin{array}{c l c l c l}
 & [{\ftwo}\op] \cdot ([\omega]\boxtimes [\omega]) & \stackrel{\ref{VER1_DEFI_SE_UND_HE}}{=} & [\he_{2,1}\op] \cdot [\omega]^{\boxtimes 2} & \stackrel{\ref{VER1_DEFI_STRUKTURMORPHISMUS_FÜR_FREEOX}}{=} & \he_{2,1}\op\mathfrak{f}_{0,X_{\operatorname{G}}} \cdot [\omega]^{\boxtimes 2}\\
 \stackrel{\text{Branching}}{=} & [\omega] \cdot \he_{2,1}\op\mathfrak{f}_{0,X_{\operatorname{G}}} & \stackrel{\ref{VER1_DEFI_STRUKTURMORPHISMUS_FÜR_FREEOX}}{=} & [\omega] \cdot [\he_{2,1}\op] & = & [\omega] \cdot [{\ftwo}\op].
\end{array}
\]

Ad 5: Note that $|\Freeopreg(1,0)|  =  1$; \cf \bfcite{Eggert}{Rem.\ 6.4}. Since $[\omega] \cdot [{\fzero}\op], [{\fzero}\op]  \in  \Freeog(1,0)$, we conclude that in $\Freeoprex$ we have
\[
\begin{array}{r c l}
[\omega] \cdot [{\fzero}\op]  & = & [{\fzero}\op].
\end{array}
\]
\epro
\elem

\blem \label{VER1_LEM_EINS_IN_GRUPPEN_EINDEUTIG}
\textit{Let $\eta'\in \Groupfp(0,1)$ such that $(\eta' \boxtimes [[\id_1\op]]) \cdot [[\mu]]  =  [[\id_1\op]]$. Then we have}
\[
\begin{array}{r c c c l}
\eta' & = & [[\eta]] & \in & \Groupfp.
\end{array}
\]
\bpro
In $\Groupfp$ we have 
\[
\begin{array}{c l c l}
 & \eta' & = & (\eta' \boxtimes [[\id_0\op]]) \cdot [[\id_1\op]]  \\
 \stackrel{\ref{VER1_DEFI_XG_UND_YG}}{=} & (\eta' \boxtimes [[\id_0\op]]) \cdot ([[\id_1\op]] \boxtimes [[\eta]]) \cdot [[\mu]] &  = & (\eta' \boxtimes [[\eta]]) \cdot [[\mu]] \\
  = & ([[\id_0\op]] \boxtimes [[\eta]]) \cdot (\eta' \boxtimes [[\id_1\op]]) \cdot [[\mu]]  & \stackrel{\nV}{=} & ([[\id_0\op]] \boxtimes [[\eta]]) \cdot [[\id_1\op]] \\
 = & [[\eta]].
\end{array}
\]
\epro
\elem

\blem \label{VER1_LEM_INVERSE_IN_GRUPPEN_EINDEUTIG}
\textit{Let $\omega'\in \Groupfp(1,1)$ such that $[[{\ftwo}\op]] \cdot (\omega' \boxtimes [[\id_1\op]]) \cdot [[\mu]]  =  [[{\fzero}\op]] \cdot [[\eta]]$. Then we have}
\[
\begin{array}{r c c c l}
\omega' & = & [[\omega]] & \in & \Groupfp.
\end{array}
\]
\bpro
In $\Groupfp$ we have
\[
\begin{array}{c l}
 & \omega' \\
 = & [[\id_1\op]] \cdot (\omega'\boxtimes [[\id_0\op]]) \\
 \stackrel{\ref{VER1_LEM_LINKSAXIOME_BEI_GRUPPEN_REICHEN}.3}{=} & [[{\ftwo}\op]] \cdot \big([[\id_1\op]] \boxtimes [[{\fzero}\op]]\big) \cdot (\omega'\boxtimes [[\id_0\op]]) \\
 = & [[{\ftwo}\op]] \cdot \big(\omega' \boxtimes [[{\fzero}\op]]\big) \\
 = & [[{\ftwo}\op]] \cdot (\omega' \boxtimes [[\id_1\op]]) \cdot \big([[\id_1\op]] \boxtimes [[{\fzero}\op]]\big) \\
 = & [[{\ftwo}\op]] \cdot (\omega' \boxtimes [[\id_1\op]]) \cdot \big([[\id_1\op]] \boxtimes [[{\fzero}\op]]\big) \cdot [[\id_1\op]] \\
 \stackrel{\ref{VER1_DEFI_XG_UND_YG}}{=} & [[{\ftwo}\op]] \cdot (\omega' \boxtimes [[\id_1\op]]) \cdot \big([[\id_1\op]] \boxtimes [[{\fzero}\op]]\big) \cdot ([[\id_1\op]] \boxtimes [[\eta]]) \cdot [[\mu]] \\
 = & [[{\ftwo}\op]] \cdot (\omega' \boxtimes [[\id_1\op]]) \cdot \Big([[\id_1\op]] \boxtimes \big([[{\fzero}\op]] \cdot [[\eta]]\big)\Big) \cdot [[\mu]] \\
 \stackrel{\ref{VER1_DEFI_XG_UND_YG}}{=} & [[{\ftwo}\op]] \cdot (\omega' \boxtimes [[\id_1\op]]) \cdot \Big([[\id_1\op]] \boxtimes \big([[{\ftwo}\op]] \cdot ([[\id_1\op]] \boxtimes [[\omega]]) \cdot [[\mu]]\big)\Big) \cdot [[\mu]] \\
 = & [[{\ftwo}\op]] \cdot (\omega' \boxtimes [[\id_1\op]]) \cdot \Big([[\id_1\op]] \boxtimes \big([[{\ftwo}\op]] \cdot ([[\id_1\op]] \boxtimes [[\omega]])\big)\Big) \cdot ([[\id_1\op]] \boxtimes [[\mu]]) \cdot [[\mu]] \\
 = & [[{\ftwo}\op]] \cdot \Big(\omega' \boxtimes \big([[{\ftwo}\op]] \cdot ([[\id_1\op]] \boxtimes [[\omega]])\big)\Big) \cdot ([[\id_1\op]] \boxtimes [[\mu]]) \cdot [[\mu]] \\
 = & [[{\ftwo}\op]] \cdot \big([[\id_1\op]] \boxtimes [[{\ftwo}\op]]\big) \cdot (\omega' \boxtimes [[\id_1\op]] \boxtimes [[\omega]]) \cdot ([[\id_1\op]] \boxtimes [[\mu]]) \cdot [[\mu]] \\
 \stackrel{\ref{VER1_LEM_LINKSAXIOME_BEI_GRUPPEN_REICHEN}.1}{=} & [[{\ftwo}\op]] \cdot \big([[{\ftwo}\op]] \boxtimes [[\id_1\op]]\big) \cdot (\omega' \boxtimes [[\id_1\op]] \boxtimes [[\omega]]) \cdot ([[\id_1\op]] \boxtimes [[\mu]]) \cdot [[\mu]] \\
 \stackrel{\ref{VER1_DEFI_XG_UND_YG}}{=} & [[{\ftwo}\op]] \cdot \big([[{\ftwo}\op]] \boxtimes [[\id_1\op]]\big) \cdot (\omega' \boxtimes [[\id_1\op]] \boxtimes [[\omega]]) \cdot ([[\mu]] \boxtimes [[\id_1\op]]) \cdot [[\mu]] \\
 = & [[{\ftwo}\op]] \cdot \Big(\big([[{\ftwo}\op]] \cdot (\omega' \boxtimes [[\id_1\op]]) \cdot [[\mu]]\big) \boxtimes [[\omega]]\Big) \cdot [[\mu]] \\
 \stackrel{\nV}{=} &  [[{\ftwo}\op]] \cdot \Big(\big([[{\fzero}\op]] \cdot [[\eta]]\big) \boxtimes [[\omega]]\Big) \cdot [[\mu]] \\
 = & [[{\ftwo}\op]] \cdot \big([[{\fzero}\op]] \boxtimes [[\id_1\op]]\big) \cdot ([[\id_0\op]]\boxtimes [[\omega]]) \cdot ([[\eta]]\boxtimes [[\id_1\op]])  \cdot [[\mu]] \\
 = & [[{\ftwo}\op]] \cdot \big([[{\fzero}\op]] \boxtimes [[\id_1\op]]\big) \cdot [[\omega]] \cdot ([[\eta]]\boxtimes [[\id_1\op]])  \cdot [[\mu]] \\
 \stackrel{\ref{VER1_LEM_LINKSAXIOME_BEI_GRUPPEN_REICHEN}.3}{=} & [[\id_1\op]] \cdot [[\omega]] \cdot ([[\eta]]\boxtimes [[\id_1\op]])  \cdot [[\mu]] \\
 = & [[\omega]] \cdot ([[\eta]]\boxtimes [[\id_1\op]])  \cdot [[\mu]] \\
 \stackrel{\ref{VER1_DEFI_XG_UND_YG}}{=} & [[\omega]] \cdot [[\id_1\op]] \\
 = & [[\omega]].
\end{array}
\]
\epro
\elem

\blem \label{VER1_LEM_INVERSE_DER_EINS_IST_EINS_IN_GRUPPEN}
\textit{We have}
\[
\begin{array}{r c c c l}
[[\eta]] \cdot [[\omega]] & = & [[\eta]] & \in & \Groupfp.
\end{array}
\]
\bpro
Note that $\id_0  =  \he_{2,0}  \in   \mapo(0,0)$ and $\he_{2,1}  =  \ftwo  \in  \mapo(2,1)$. Note that \mbox{$|\Groupfp(0,0)|  =  1$;} \cf \bfcite{Eggert}{Rem.\ 6.4}.
Since $[[\eta]] \cdot [[{\fzero}\op]], [[\id_0\op]]  \in  \Groupfp(0,0)$, we conclude that
\[
\begin{array}{r c c c l}
[[\eta]] \cdot [[{\fzero}\op]] & = &  [[\id_0\op]] & \in & \Groupfp(0,0).
\end{array}
\]
In $\Groupfp$ we have
\[
\begin{array}{c l c l}
 & [[\eta]] \cdot [[\omega]] & = & \big(([[\eta]] \cdot [[\omega]]) \boxtimes [[\id_0\op]]\big) \cdot [[\id_1\op]] \\
 \stackrel{\ref{VER1_DEFI_XG_UND_YG}}{=} & \big(([[\eta]] \cdot [[\omega]]) \boxtimes [[\id_0\op]]\big) \cdot ([[\id_1\op]] \boxtimes [[\eta]]) \cdot [[\mu]]  & = & \big(([[\eta]] \cdot [[\omega]]) \boxtimes [[\eta]]\big) \cdot [[\mu]] \\
 = & ([[\eta]] \boxtimes [[\eta]]) \cdot ([[\omega]] \boxtimes [[\id_1\op]]) \cdot [[\mu]]  & = & [[\id_0\op]] \cdot ([[\eta]] \boxtimes [[\eta]]) \cdot ([[\omega]] \boxtimes [[\id_1\op]]) \cdot [[\mu]] \\
 = & [[\he_{2,0}\op]] \cdot [[\eta]]^{\boxtimes 2} \cdot ([[\omega]] \boxtimes [[\id_1\op]]) \cdot [[\mu]]  & = & \he_{2,0}\op\mathfrak{g}_0\fp \cdot [[\eta]]^{\boxtimes 2} \cdot ([[\omega]] \boxtimes [[\id_1\op]]) \cdot [[\mu]] \\
 \stackrel{\text{Braiding}}{=} & [[\eta]] \cdot \he_{2,1}\op\mathfrak{g}_0\fp \cdot ([[\omega]] \boxtimes [[\id_1\op]]) \cdot [[\mu]] & = & [[\eta]] \cdot [[\he_{2,1}\op]] \cdot ([[\omega]] \boxtimes [[\id_1\op]]) \cdot [[\mu]] \\
 = & [[\eta]] \cdot [[{\ftwo}\op]] \cdot ([[\omega]] \boxtimes [[\id_1\op]]) \cdot [[\mu]]  & \stackrel{\ref{VER1_DEFI_XG_UND_YG}}{=}& [[\eta]] \cdot [[{\fzero}\op]] \cdot [[\eta]] \\
 \stackrel{\text{see above}}{=} & [[\id_0\op]] \cdot [[\eta]]  & = & [[\eta]].
\end{array}
\]
\epro
\elem

\blem \label{VER1_LEM_INVERSE_DER_INVERSEN_IST_EINS_IN_GRUPPEN}
\textit{We have}
\[
\begin{array}{r c c c l}
[[\omega]] \cdot [[\omega]] & = & [[\id_1\op]] & \in & \Groupfp.
\end{array}
\]
\bpro
In $\Groupfp$ we have
\[
\begin{array}{c l}
 & [[\omega]] \cdot [[\omega]] \\
 = & [[\id_1\op]] \cdot \big(([[\omega]] \cdot [[\omega]])\boxtimes [[\id_0\op]]\big) \cdot [[\id_1\op]] \\
 \stackrel{\ref{VER1_LEM_LINKSAXIOME_BEI_GRUPPEN_REICHEN}.3}{=} & [[{\ftwo}\op]] \cdot \big([[\id_1\op]] \boxtimes [[{\fzero}\op]]\big) \cdot \big(([[\omega]] \cdot [[\omega]])\boxtimes [[\id_0\op]]\big) \cdot [[\id_1\op]] \\
 \stackrel{\ref{VER1_DEFI_XG_UND_YG}}{=} & [[{\ftwo}\op]] \cdot \big([[\id_1\op]] \boxtimes [[{\fzero}\op]]\big) \cdot \big(([[\omega]] \cdot [[\omega]])\boxtimes [[\id_0\op]]\big) \cdot ([[\id_1\op]] \boxtimes [[\eta]]) \cdot [[\mu]] \\
 = & [[{\ftwo}\op]] \cdot \Big(([[\omega]] \cdot [[\omega]])\boxtimes \big([[{\fzero}\op]] \cdot [[\eta]]\big)\Big) \cdot [[\mu]] \\
 \stackrel{\ref{VER1_DEFI_XG_UND_YG}}{=} & [[{\ftwo}\op]] \cdot \Big(([[\omega]] \cdot [[\omega]])\boxtimes \big([[{\ftwo}\op]] \cdot ([[\omega]] \boxtimes [[\id_1\op]]) \cdot [[\mu]]\big)\Big) \cdot [[\mu]] \\
 = & [[{\ftwo}\op]] \cdot \big([[\id_1\op]] \boxtimes [[{\ftwo}\op]]\big) \cdot \big(([[\omega]] \cdot [[\omega]]) \boxtimes [[\omega]] \boxtimes [[\id_1\op]]\big) \cdot  ([[\id_1\op]] \boxtimes [[\mu]]) \cdot [[\mu]] \\
 \stackrel{\ref{VER1_LEM_LINKSAXIOME_BEI_GRUPPEN_REICHEN}.1}{=} & [[{\ftwo}\op]] \cdot \big([[{\ftwo}\op]] \boxtimes [[\id_1\op]]\big) \cdot \big(([[\omega]] \cdot [[\omega]]) \boxtimes [[\omega]] \boxtimes [[\id_1\op]]\big) \cdot  ([[\id_1\op]] \boxtimes [[\mu]]) \cdot [[\mu]] \\
 \stackrel{\ref{VER1_DEFI_XG_UND_YG}}{=} & [[{\ftwo}\op]] \cdot \big([[{\ftwo}\op]] \boxtimes [[\id_1\op]]\big) \cdot \big(([[\omega]] \cdot [[\omega]]) \boxtimes [[\omega]] \boxtimes [[\id_1\op]]\big) \cdot  ([[\mu]] \boxtimes [[\id_1\op]]) \cdot [[\mu]] \\
 = & [[{\ftwo}\op]] \cdot \big([[{\ftwo}\op]] \boxtimes [[\id_1\op]]\big) \cdot ([[\omega]]\boxtimes [[\omega]] \boxtimes [[\id_1\op]]) \cdot ([[\omega]] \boxtimes [[\id_1\op]] \boxtimes [[\id_1\op]]) \\
 & \cdot \, ([[\mu]] \boxtimes [[\id_1\op]]) \cdot [[\mu]] \\
 = & [[{\ftwo}\op]] \cdot \Big(\big([[{\ftwo}\op]] \cdot ([[\omega]]\boxtimes [[\omega]])\big) \boxtimes [[\id_1\op]]\Big) \cdot ([[\omega]] \boxtimes [[\id_1\op]] \boxtimes [[\id_1\op]])  \\
 & \cdot \, ([[\mu]] \boxtimes [[\id_1\op]]) \cdot [[\mu]] \\
  \stackrel{\ref{VER1_LEM_LINKSAXIOME_BEI_GRUPPEN_REICHEN}.4}{=} & [[{\ftwo}\op]] \cdot \Big(\big([[\omega]] \cdot [[{\ftwo}\op]]\big) \boxtimes [[\id_1\op]]\Big) \cdot ([[\omega]] \boxtimes [[\id_1\op]] \boxtimes [[\id_1\op]]) \\
  &  \cdot \, ([[\mu]] \boxtimes [[\id_1\op]]) \cdot [[\mu]] \\
 = & [[{\ftwo}\op]] \cdot ([[\omega]]\boxtimes [[\id_1\op]]) \cdot \big([[{\ftwo}\op]] \boxtimes [[\id_1\op]]\big) \cdot ([[\omega]] \boxtimes [[\id_1\op]] \boxtimes [[\id_1\op]]) \\
 &  \cdot \, ([[\mu]] \boxtimes [[\id_1\op]]) \cdot [[\mu]] 
 \end{array}
\]
\[
\begin{array}{c l}
\hspace{-6.3ex} = & [[{\ftwo}\op]] \cdot ([[\omega]]\boxtimes [[\id_1\op]]) \cdot \Big(\big([[{\ftwo}\op]] \cdot ([[\omega]] \boxtimes [[\id_1\op]]) \cdot [[\mu]]\big) \boxtimes [[\id_1\op]]\Big) \cdot [[\mu]] \\
\hspace{-6.3ex} \stackrel{\ref{VER1_DEFI_XG_UND_YG}}{=} & [[{\ftwo}\op]] \cdot ([[\omega]]\boxtimes [[\id_1\op]]) \cdot \Big(\big([[{\fzero}\op]] \cdot [[\eta]]\big) \boxtimes [[\id_1\op]]\Big) \cdot [[\mu]] \\
\hspace{-6.3ex} = & [[{\ftwo}\op]] \cdot ([[\omega]]\boxtimes [[\id_1\op]]) \cdot \big([[{\fzero}\op]] \boxtimes [[\id_1\op]]\big) \cdot ([[\eta]] \boxtimes [[\id_1\op]]) \cdot [[\mu]] \\
\hspace{-6.3ex} = & [[{\ftwo}\op]] \cdot \Big(\big([[\omega]] \cdot [[{\fzero}\op]]\big) \boxtimes [[\id_1\op]]\Big) \cdot ([[\eta]] \boxtimes [[\id_1\op]]) \cdot [[\mu]] \\
\hspace{-6.3ex} \stackrel{\ref{VER1_LEM_LINKSAXIOME_BEI_GRUPPEN_REICHEN}.5}{=} & [[{\ftwo}\op]] \cdot \big([[{\fzero}\op]] \boxtimes [[\id_1\op]]\big) \cdot ([[\eta]] \boxtimes [[\id_1\op]]) \cdot [[\mu]] \\
\hspace{-6.3ex} \stackrel{\ref{VER1_LEM_LINKSAXIOME_BEI_GRUPPEN_REICHEN}.3}{=} & [[\id_1\op]] \cdot ([[\eta]] \boxtimes [[\id_1\op]]) \cdot [[\mu]] \\
\hspace{-6.3ex} = & ([[\eta]] \boxtimes [[\id_1\op]]) \cdot [[\mu]] \\
\hspace{-6.3ex} \stackrel{\ref{VER1_DEFI_XG_UND_YG}}{=} & [[\id_1\op]].
\end{array}
\]
\epro
\elem

\bbem
Lemma \ref{VER1_LEM_EINS_IN_GRUPPEN_EINDEUTIG} is the operad version of the fact that the neutral element of a group is unique.

Lemma \ref{VER1_LEM_INVERSE_IN_GRUPPEN_EINDEUTIG} is the operad version of the fact that the inverse of a group element is unique.

Lemma \ref{VER1_LEM_INVERSE_DER_EINS_IST_EINS_IN_GRUPPEN} is the operad version of the fact that the inverse of the neutral element of a group is again the neutral element.

Lemma \ref{VER1_LEM_INVERSE_DER_INVERSEN_IST_EINS_IN_GRUPPEN} is the operad version of the fact that the inverse of the inverse of a group element
is again the group element.

In each case it was possible to translate the usual proof into an operadic version.
\ebem

\bdefi
We define the subset 
\[
\begin{array}{r c c c c c l}
Z_{\operatorname{G}} & := & \{(y_{\operatorname{G},1}^1,y_{\operatorname{G},2}^1), (y_{\operatorname{G},1}^2,y_{\operatorname{G},2}^2), (y_{\operatorname{G},1}^4,y_{\operatorname{G},2}^4)\} & \subseteq & Y_{\operatorname{G}} & \subseteq & \Freeopreg \btimes \Freeopreg;
\end{array}
\]
\cf definition \ref{VER1_DEFI_BIINDEXED_SETS}, \cf definition \ref{VER1_DEFI_XG_UND_YG}.
\edefi

\blem \label{VER1_LEM_YG_GLEICH_ZG}
\textit{We have}
\[
\begin{array}{r c c c l}
(\equiv_{Z_{\operatorname{G}}}) & = & (\equiv_{Y_{\operatorname{G}}}) & \subseteq & \Freeopreg \btimes \Freeopreg;
\end{array}
\]
\cf definition \ref{VER1_DEFI_EQUIV_Y}.
\bpro
Since $Z_{\operatorname{G}}  \subseteq  Y_{\operatorname{G}}$, we conclude that $(\equiv_{Z_{\operatorname{G}}})  \subseteq  (\equiv_{Y_{\operatorname{G}}})$; \cf definition \ref{VER1_DEFI_EQUIV_Y}. To prove that $(\equiv_{Y_{\operatorname{G}}})  \subseteq  (\equiv_{Z_{\operatorname{G}}})$ it suffices to show that $(y_{\operatorname{G},1}^3,y_{\operatorname{G},2}^3), (y_{\operatorname{G},1}^5,y_{\operatorname{G},2}^5)  \in  (\equiv_{Z_{\operatorname{G}}})$.\\

\textit{Claim 1: We have $(y_{\operatorname{G},1}^5,y_{\operatorname{G},2}^5) \in (\equiv_{Z_{\operatorname{G}}})$, $\ie$} 
\[
\begin{array}{r c c c c c l}
y_{\operatorname{G},1}^5 & = & [{\ftwo}\op] \cdot ([\id_1\op]\boxtimes [\omega]) \cdot [\mu] & \equiv_{Z_{\operatorname{G}}} & [{\fzero}\op]  \cdot [\eta] & = & y_{\operatorname{G},2}^5.
\end{array}
\]

In $\Freeopreg$ we have 
\[
\begin{array}{c l}
 & y_{\operatorname{G},1}^5 \\
 = & [{\ftwo}\op] \cdot ([\id_1\op]\boxtimes [\omega]) \cdot [\mu] \\
 = & [\id_1\op] \cdot [{\ftwo}\op] \cdot ([\id_1\op]\boxtimes [\omega]) \cdot [\mu] \\
 \stackrel{\ref{VER1_LEM_LINKSAXIOME_BEI_GRUPPEN_REICHEN}.3}{=} & [{\ftwo}\op] \cdot \big([{\fzero}\op] \boxtimes [\id_1\op]\big)  \cdot [{\ftwo}\op] \cdot ([\id_1\op]\boxtimes [\omega]) \cdot [\mu] \\
 = & [{\ftwo}\op] \cdot \big([{\fzero}\op] \boxtimes [\id_1\op]\big)  \cdot \big([\id_0\op] \boxtimes [{\ftwo}\op]\big) \cdot ([\id_1\op]\boxtimes [\omega]) \cdot [\mu] \\
 = & [{\ftwo}\op] \cdot  \big([{\fzero}\op] \boxtimes [{\ftwo}\op]\big) \cdot ([\id_1\op]\boxtimes [\omega]) \cdot [\mu] \\
 \hspace{2.6ex} \equiv_{Z_{\operatorname{G}}} & [{\ftwo}\op] \cdot  \big([{\fzero}\op] \boxtimes [{\ftwo}\op]\big) \cdot \Big(\big(([\eta] \boxtimes [\id_1\op]) \cdot [\mu]\big)\boxtimes [\omega]\Big) \cdot [\mu] \\
 = & [{\ftwo}\op] \cdot  \big([{\fzero}\op] \boxtimes [{\ftwo}\op]\big) \cdot \Big(\big(([\eta] \boxtimes [\id_1\op]) \cdot [\mu]\big)\boxtimes ([\omega] \cdot [\id_1\op])\Big) \cdot [\mu] \\
 = & [{\ftwo}\op] \cdot \big([{\fzero}\op] \boxtimes [{\ftwo}\op]\big) \cdot ([\eta] \boxtimes [\id_1\op] \boxtimes [\omega]) \cdot ([\mu]\boxtimes [\id_1\op]) \cdot [\mu] \\
 = & [{\ftwo}\op] \cdot \Big(\big([{\fzero}\op] \cdot [\eta]\big)  \boxtimes  \big([{\ftwo}\op] \cdot ([\id_1\op] \boxtimes [\omega])\big)\Big) \cdot ([\mu]\boxtimes [\id_1\op]) \cdot [\mu] \\
 \stackrel{\ref{VER1_LEM_LINKSAXIOME_BEI_GRUPPEN_REICHEN}.5}{=} & [{\ftwo}\op] \cdot \Big(\big([\omega] \cdot [{\fzero}\op] \cdot [\eta]\big)  \boxtimes  \big([{\ftwo}\op] \cdot ([\id_1\op] \boxtimes [\omega])\big)\Big) \cdot ([\mu]\boxtimes [\id_1\op]) \cdot [\mu] \\
 \hspace{2.6ex} \equiv_{Z_{\operatorname{G}}} &  [{\ftwo}\op] \cdot \Big(\big([\omega] \cdot [{\ftwo}\op] \cdot ([\omega] \boxtimes [\id_1\op]) \cdot [\mu]\big)  \boxtimes  \big([{\ftwo}\op] \cdot ([\id_1\op] \boxtimes [\omega])\big)\Big) \cdot ([\mu]\boxtimes [\id_1\op]) \cdot [\mu] \\
 \stackrel{\ref{VER1_LEM_LINKSAXIOME_BEI_GRUPPEN_REICHEN}.4}{=} & [{\ftwo}\op] \cdot \Big(\big( [{\ftwo}\op] \cdot ([\omega] \boxtimes [\omega]) \cdot ([\omega] \boxtimes [\id_1\op]) \cdot [\mu]\big)  \boxtimes  \big([{\ftwo}\op] \cdot ([\id_1\op] \boxtimes [\omega])\big)\Big) \cdot ([\mu]\boxtimes [\id_1\op]) \cdot [\mu] \\
 = & [{\ftwo}\op] \cdot \bigg(\Big([{\ftwo}\op] \cdot \big(([\omega] \cdot [\omega]) \boxtimes [\omega]\big) \cdot [\mu]\Big)  \boxtimes  \big([{\ftwo}\op] \cdot ([\id_1\op] \boxtimes [\omega]) \cdot [\id_2\op]\big)\bigg) \cdot ([\mu]\boxtimes [\id_1\op]) \cdot [\mu] \\
 = & [{\ftwo}\op] \cdot \Big(\big([{\ftwo}\op] \boxtimes [{\ftwo}\op]\big) \cdot \big(([\omega] \cdot [\omega]) \boxtimes [\omega] \boxtimes [\id_1\op] \boxtimes [\omega]\big) \cdot ([\mu]\boxtimes [\id_2\op])\Big) \cdot ([\mu]\boxtimes [\id_1\op]) \cdot [\mu] \\
 = & [{\ftwo}\op] \cdot \big([{\ftwo}\op] \boxtimes [{\ftwo}\op]\big) \cdot \big(([\omega] \cdot [\omega]) \boxtimes [\omega] \boxtimes [\id_1\op] \boxtimes [\omega]\big) \cdot \Big(\big(([\mu]\boxtimes [\id_1\op]) \cdot [\mu]\big)\boxtimes [\id_1\op]\Big) \cdot [\mu] \\
 \hspace{2.6ex} \equiv_{Z_{\operatorname{G}}} & [{\ftwo}\op] \cdot \big([{\ftwo}\op] \boxtimes [{\ftwo}\op]\big) \cdot \big(([\omega] \cdot [\omega]) \boxtimes [\omega] \boxtimes [\id_1\op] \boxtimes [\omega]\big) \cdot \Big(\big(([\id_1\op]\boxtimes [\mu]) \cdot [\mu]\big)\boxtimes [\id_1\op]\Big) \cdot [\mu] \\
 \stackrel{\ref{VER1_LEM_LINKSAXIOME_BEI_GRUPPEN_REICHEN}.2}{=} & [{\ftwo}\op] \cdot \big([{\ftwo}\op] \boxtimes [\id_1\op]\big) \cdot \big([\id_1\op] \boxtimes [{\ftwo}\op] \boxtimes [\id_1\op]\big) \cdot \big(([\omega] \cdot [\omega]) \boxtimes [\omega] \boxtimes [\id_1\op] \boxtimes [\omega]\big) \\
  & \cdot \, \Big(\big(([\id_1\op]\boxtimes [\mu]) \cdot [\mu]\big)\boxtimes [\id_1\op]\Big) \cdot [\mu] \\
 = & [{\ftwo}\op] \cdot \big([{\ftwo}\op] \boxtimes [\id_1\op]\big) \cdot \Big(([\omega] \cdot [\omega]) \boxtimes \big([{\ftwo}\op] \cdot ([\omega] \boxtimes [\id_1\op])\big) \boxtimes [\omega]\Big) \\
  &  \cdot \, ([\id_1\op]\boxtimes [\mu] \boxtimes [\id_1\op]) \cdot ([\mu] \boxtimes [\id_1\op]) \cdot [\mu] \\
 = & [{\ftwo}\op] \cdot \big([{\ftwo}\op] \boxtimes [\id_1\op]\big) \cdot \Big(([\omega] \cdot [\omega]) \boxtimes \big([{\ftwo}\op] \cdot ([\omega] \boxtimes [\id_1\op]) \cdot [\mu]\big) \boxtimes [\omega]\Big) \cdot ([\mu] \boxtimes [\id_1\op]) \cdot [\mu] \\
 \hspace{2.6ex} \equiv_{Z_{\operatorname{G}}} & [{\ftwo}\op] \cdot \big([{\ftwo}\op] \boxtimes [\id_1\op]\big) \cdot \Big(([\omega] \cdot [\omega]) \boxtimes \big([{\fzero}\op] \cdot [\eta]\big) \boxtimes [\omega]\Big) \cdot ([\mu] \boxtimes [\id_1\op]) \cdot [\mu] \\
 \hspace{2.6ex} \equiv_{Z_{\operatorname{G}}} & [{\ftwo}\op] \cdot \big([{\ftwo}\op] \boxtimes [\id_1\op]\big) \cdot \Big(([\omega] \cdot [\omega]) \boxtimes \big([{\fzero}\op] \cdot [\eta]\big) \boxtimes [\omega]\Big) \cdot ([\id_1\op] \boxtimes [\mu]) \cdot [\mu] \\
 = & [{\ftwo}\op] \cdot \big([{\ftwo}\op] \boxtimes [\id_1\op]\big) \cdot \big(([\omega] \cdot [\omega]) \boxtimes [{\fzero}\op] \boxtimes [\omega]\big) \cdot ([\id_1\op] \boxtimes [\eta] \boxtimes [\id_1\op]) \cdot ([\id_1\op] \boxtimes [\mu]) \cdot [\mu] \\
 = & [{\ftwo}\op] \cdot \big([{\ftwo}\op] \boxtimes [\id_1\op]\big) \cdot \big(([\omega] \cdot [\omega]) \boxtimes [{\fzero}\op] \boxtimes [\omega]\big) \cdot \Big([\id_1\op] \boxtimes \big(([\eta] \boxtimes [\id_1\op]) \cdot [\mu]\big)\Big) \cdot [\mu] \\
 \hspace{2.6ex} \equiv_{Z_{\operatorname{G}}} & [{\ftwo}\op] \cdot  \big([{\ftwo}\op] \boxtimes [\id_1\op]\big) \cdot \big(([\omega] \cdot [\omega]) \boxtimes [{\fzero}\op] \boxtimes [\omega]\big) \cdot ([\id_1\op] \boxtimes [\id_1\op]) \cdot [\mu] \\
 = & [{\ftwo}\op] \cdot  \big([{\ftwo}\op] \boxtimes [\id_1\op]\big) \cdot \big(([\omega] \cdot [\omega]) \boxtimes [{\fzero}\op] \boxtimes [\omega]\big) \cdot [\mu] \\
 = & [{\ftwo}\op] \cdot \bigg(\Big([{\ftwo}\op] \cdot \big(([\omega] \cdot [\omega]) \boxtimes [{\fzero}\op]\big)\Big) \boxtimes [\omega]\bigg)  \cdot [\mu] \\
 = & [{\ftwo}\op] \cdot \bigg(\Big([{\ftwo}\op] \cdot \big([\id_1\op] \boxtimes [{\fzero}\op]\big) \cdot \big(([\omega] \cdot [\omega]) \boxtimes [\id_0\op]\big)\Big) \boxtimes [\omega]\bigg)  \cdot [\mu] \\ 
 = & [{\ftwo}\op] \cdot \bigg(\Big([{\ftwo}\op] \cdot \big([\id_1\op] \boxtimes [{\fzero}\op]\big) \cdot [\omega] \cdot [\omega]\Big) \boxtimes [\omega]\bigg)  \cdot [\mu] 
 \end{array}
\]
\[
\begin{array}{c l}
 \hspace{-50ex} \stackrel{\ref{VER1_LEM_LINKSAXIOME_BEI_GRUPPEN_REICHEN}.3}{=} &   [{\ftwo}\op] \cdot \big(([\id_1\op] \cdot [\omega] \cdot [\omega]) \boxtimes [\omega]\big)  \cdot [\mu] \\
 \hspace{-50ex} = &   [{\ftwo}\op] \cdot \big(([\omega] \cdot [\omega]) \boxtimes [\omega]\big)  \cdot [\mu] \\
 \hspace{-50ex} = &   [{\ftwo}\op] \cdot ([\omega] \boxtimes [\omega]) \cdot ([\omega] \boxtimes [\id_1\op])  \cdot [\mu] \\
 \hspace{-50ex} \stackrel{\ref{VER1_LEM_LINKSAXIOME_BEI_GRUPPEN_REICHEN}.4}{=} &   [\omega] \cdot [{\ftwo}\op]  \cdot ([\omega] \boxtimes [\id_1\op])  \cdot [\mu] \\
 \hspace{-50ex} \hspace{2.6ex} \equiv_{Z_{\operatorname{G}}} &   [\omega] \cdot [{\fzero}\op]  \cdot [\eta] \\
 \hspace{-50ex} \stackrel{\ref{VER1_LEM_LINKSAXIOME_BEI_GRUPPEN_REICHEN}.5}{=} &   [{\fzero}\op]  \cdot [\eta]\\
 \hspace{-50ex} = &   y_{\operatorname{G},2}^5.
\end{array}
\]
This proves \textit{Claim 1}.\\

\textit{Claim 2: We have $(y_{\operatorname{G},1}^3,y_{\operatorname{G},2}^3) \in (\equiv_{Z_{\operatorname{G}}})$, $\ie$}
\[
\begin{array}{r c c c c c l}
y_{\operatorname{G},1}^3 & = & ([\id_1] \boxtimes [\eta]) \cdot [\mu] & \equiv_{Z_{\operatorname{G}}} & [\id_1\op] & = & y_{\operatorname{G},2}^3.
\end{array}
\]

By Claim 1 in $\Freeog$ we have
\[
\begin{array}{r c l}
y_{\operatorname{G},1}^3  & = & ([\id_1\op] \boxtimes [\eta])\cdot [\mu] \\
 & = & [\id_1\op] \cdot ([\id_1\op] \boxtimes [\eta])\cdot [\mu] \\
 & \stackrel{\ref{VER1_LEM_LINKSAXIOME_BEI_GRUPPEN_REICHEN}.3}{=} & [{\ftwo}\op] \cdot \big([\id_1\op] \boxtimes [{\fzero}\op]\big) \cdot ([\id_1\op] \boxtimes [\eta])\cdot [\mu] \\
 & = & [{\ftwo}\op] \cdot \Big([\id_1\op] \boxtimes \big([{\fzero}\op] \cdot [\eta]\big)\Big) \cdot  [\mu] \\
 & \hspace{2.6ex} \equiv_{Z_{\operatorname{G}}} & [{\ftwo}\op] \cdot \Big([\id_1\op] \boxtimes \big([{\ftwo}\op] \cdot ([\omega] \boxtimes [\id_1\op]) \cdot [\mu]\big)\Big) \cdot  [\mu] \\
 & = & [{\ftwo}\op] \cdot \big([\id_1\op] \boxtimes [{\ftwo}\op]\big) \cdot ([\id_1\op] \boxtimes [\omega] \boxtimes [\id_1\op]) \cdot ([\id_1\op] \boxtimes [\mu]) \cdot [\mu] \\
 & \stackrel{\ref{VER1_LEM_LINKSAXIOME_BEI_GRUPPEN_REICHEN}.1}{=} & [{\ftwo}\op] \cdot \big([{\ftwo}\op] \boxtimes [\id_1\op]\big) \cdot ([\id_1\op] \boxtimes [\omega] \boxtimes [\id_1\op]) \cdot ([\id_1\op] \boxtimes [\mu]) \cdot [\mu] \\
 & \hspace{2.6ex} \equiv_{Z_{\operatorname{G}}} & [{\ftwo}\op] \cdot \big([{\ftwo}\op] \boxtimes [\id_1\op]\big) \cdot ([\id_1\op] \boxtimes [\omega] \boxtimes [\id_1\op]) \cdot ([\mu] \boxtimes [\id_1\op]) \cdot [\mu] \\
 & = & [{\ftwo}\op] \cdot \Big(\big([{\ftwo}\op] \cdot ([\id_1\op] \boxtimes [\omega]) \cdot [\mu]\big) \boxtimes [\id_1\op]\Big)  \cdot [\mu] \\
 & \hspace{2.6ex} \stackrel{\hspace{-2.5ex}\text{Claim 1}}{\equiv_{Z_{\operatorname{G}}}} & [{\ftwo}\op] \cdot \Big(\big([{\fzero}\op] \cdot [\eta]\big) \boxtimes [\id_1\op]\Big)  \cdot [\mu] \\
 & = & [{\ftwo}\op] \cdot \big([{\fzero}\op] \boxtimes [\id_1\op]\big) \cdot ([\eta] \boxtimes [\id_1\op]) \cdot [\mu] \\
 & \stackrel{\ref{VER1_LEM_LINKSAXIOME_BEI_GRUPPEN_REICHEN}.3}{=} & [\id_1\op] \cdot ([\eta] \boxtimes [\id_1\op]) \cdot [\mu] \\
 & = & ([\eta] \boxtimes [\id_1\op]) \cdot [\mu] \\
 & \hspace{2.6ex} \equiv_{Z_{\operatorname{G}}} & [\id_1\op] \\
 & = & y_{\operatorname{G},2}^3.
\end{array}
\]
This proves \textit{Claim 2}.
\epro
\elem

\bpropo \label{VER1_PROPO_GROUPFP_GEHT_AUCH_MIT_DREI_AXIOMEN}
\textit{We have}
\[
\begin{array}{r c c c l}
\Groupfp & = & \generated\xgyg & = & \generated\big<X_{\operatorname{G}}\,|\,Z_{\operatorname{G}}\big>.
\end{array}
\]
\bpro
We have
\[
\begin{array}{r c c c c c c c l}
\Groupfp \hspace{-0.4ex} & \hspace{-0.4ex} \stackrel{\text{\ref{VER1_DEFI_XG_UND_YG}}}{=} \hspace{-0.4ex} & \hspace{-0.4ex} \generated\xgyg \hspace{-0.4ex} & \hspace{-0.4ex} \stackrel{\text{\ref{VER1_DEFI_OPERADE_X_MODULO_Y}}}{=} \hspace{-0.4ex} & \hspace{-0.4ex} \Freeog/(\equiv_{Y_{\operatorname{G}}}) \hspace{-0.4ex} & \hspace{-0.4ex} \stackrel{\ref{VER1_LEM_YG_GLEICH_ZG}}{=} \hspace{-0.4ex} & \hspace{-0.4ex} \Freeog/(\equiv_{Z_{\operatorname{G}}}) \hspace{-0.4ex} & \hspace{-0.4ex} \stackrel{\text{\ref{VER1_DEFI_OPERADE_X_MODULO_Y}}}{=} \hspace{-0.4ex} & \hspace{-0.4ex} \generated\big<X_{\operatorname{G}}\,|\,Z_{\operatorname{G}}\big>.
\end{array}
\]
\epro
\epropo

\bbem
Proposition \ref{VER1_PROPO_GROUPFP_GEHT_AUCH_MIT_DREI_AXIOMEN} is the operad version of the fact that the associativity of the multiplication, the left-neutrality of the unit and the existence of left-inverse elements are sufficient to define a group in the usual sense; \cfeg \bfcite{Kurosh}{p.\ 33--34}.

It was possible to translate the usual proof into an operadic version.
\ebem

\begin{footnotesize}

\end{footnotesize}


\begin{thebibliography}{999}

%\bibitem[AL]{AguiarLivernet} Aguiar, M.; Livernet, M., \textit{The associative operad and the weak order on the symmetric groups}, Arxiv math/0511698v3, 2005. \url{https://arxiv.org/abs/math/0511698v3}

\bibitem[BV]{BV} Boardman, J.M.; Vogt, R.M., \textit{Homotopy-everything H-spaces}, Bull.\ Amer.\ Math.\ Soc.\ 74, S.\ 1117-1122, 1968. \url{https://doi.org/10.1090/S0002-9904-1968-12070-1}
 
\bibitem[Egg]{Eggert} Eggert, V., \textit{Operads in the sense of Mac\,Lane}, Universität Stuttgart, 2017.\\ \url{http://dx.doi.org/10.18419/opus-13251}

%\bibitem[LV]{LodayVallette} Loday, J.L.; Vallette, B., \textit{Algebraic Operads}, Grundlehren der mathematischen Wissenschaften, Volume 346, Springer-Verlag, 2012.

%\bibitem[Mar]{Markl} Markl, M., \textit{Homotopy Algebras are Homotopy Algebras}, Arxiv math/9907138v3, 1999. \\ \url{https://arxiv.org/abs/math/9907138}

\bibitem[May]{May} May, P., \textit{The geometry of iterated loop spaces}, Lecture Notes in Mathematics, Vol. 271, Springer (1972)\\
\url{http://www.math.uchicago.edu/~may/BOOKSMaster.html}

\bibitem[May2]{May2} May, P., \textit{Operads, algebras and modules}, Conference Proc.\ ``Operads: Proceedings of Renaissance Conferences'', ed.\ J.L.\ Loday, J.D.\ Stasheff, A.A.\ Voronov, S.\ 15--31,
Contemp.\ Math.\ 202, 1997.\\ \url{www.math.uchicago.edu/~may/PAPERS/mayi.pdf}

\bibitem[ML]{MacLane} Mac\,Lane, S., \textit{Categorical algebra}, Bull.\ Amer.\ Math.\ Soc.\ 71, S.\ 40--106, 1965.\\ \url{https://www.ams.org/journals/bull/1965-71-01/S0002-9904-1965-11234-4/S0002-9904-1965-11234-4.pdf}

\bibitem[ML2]{MacLanePACTs} Mac\,Lane, S., \textit{Higher homotopies, $\mathrm{PACT}$s, and the Bar Construction}, Conference Proc.\ ``Homotopy Invariant Algebraic Structures: A Conference in Honor of J.\ Michael Boardman'', ed.\ J.P.\ Meyer, J.\ Morava, W.S.\ Wilson, 
Contemp.\ Math.\ 239, S.\ 5--7, 1999. \url{http://dx.doi.org/10.1090/conm/239/3591}

%\bibitem[MP]{MorrisonPenneys} Morrison, S.; Penneys, D., \textit{Monoidal categories enriched in braided monoidal categories}, arXiv:1701.00567v1, 2017. \\ \url{https://arxiv.org/pdf/1701.00567.pdf}

\bibitem[MSS]{MarklShniderStasheff} Markl, M.; Shnider, S.; Stasheff, J.,  \textit{Operads in Algebra, Topology and Physics}, AMS Math.\ Surveys and Monographs 96, 2002.

%\bibitem[Sta2]{Sta2} Stasheff, J., \textit{What Is\,\dots\,an Operad?}, Notices AMS 51(6), S.\ 630--631.

\bibitem[Kurosh]{Kurosh} Kurosh, A.G., \textit{The theory of groups}, Chelsea Publishing Co., New York, N.Y., Vol. I Translated and edited by K.A. Hirsch, 1955

\end{thebibliography}
\end{document}